\documentstyle[amscd,amssymb,verbatim,12pt]{amsart}
\newcommand{\Aff}{\mbox{\rm Aff}}

\newcommand{\Aut}{\mbox{\rm Aut}}
\newcommand{\bb}{\mbox{\rm b}}

\newcommand{\C}{{\Bbb C}}

\newcommand{\cone}{\mbox{\rm cone}}

\newcommand{\cyl}{\mbox{\rm cyl}}

\newcommand{\diam}{\mbox{\rm diam}}
\newcommand{\Diff}{\mbox{\rm Diff}}

\newcommand{\Dom}{\mbox{\rm Dom}}
\newcommand{\End}{\mbox{\rm End}}

\newcommand{\GL}{\mbox{\rm GL}}

\newcommand{\HH}{\mbox{\rm H}}

\newcommand{\Hom}{\mbox{\rm Hom}}

\newcommand{\Id}{\mbox{\rm Id}}

\newcommand{\Image}{\mbox{\rm Im}}

\newcommand{\Int}{\mbox{\rm int}}

\newcommand{\Isom}{\mbox{\rm Isom}}
\newcommand{\Ker}{\mbox{\rm Ker}}

\newcommand{\R}{{\Bbb R}}

\newcommand{\Ric}{\mbox{\rm Ric}}

\newcommand{\SO}{\mbox{\rm SO}}
\newcommand{\SU}{\mbox{\rm SU}}

\newcommand{\U}{\mbox{\rm U}}

\newcommand{\vol}{\mbox{\rm vol}}
\newcommand{\Z}{{\Bbb Z}}
\theoremstyle{plain}
\newtheorem{definition}{Definition}

\newtheorem{lemma}{Lemma}
\newtheorem{theorem}{Theorem}
\newtheorem{proposition}{Proposition}
\newtheorem{corollary}{Corollary}

\numberwithin{equation}{section}
\renewcommand{\rm}{\normalshape}
\begin{document}
\title{Collapsing and the Differential Form Laplacian : The Case of a
Singular Limit Space}
\author{John Lott}
\address{Department of Mathematics\\
University of Michigan\\
Ann Arbor, MI  48109-1109\\
USA}
\email{lott@@math.lsa.umich.edu}
\thanks{Research supported by NSF grant DMS-0072154}
\subjclass{Primary: 58G25; Secondary: 53C23}
\date{January 28, 2002}
\maketitle
\begin{abstract}
We analyze the limit of the $p$-form
Laplacian under a collapse with bounded sectional curvature and bounded
diameter to a singular limit space.  
As applications, we give results about upper and
lower bounds on the $j$-th eigenvalue
of the $p$-form Laplacian, in terms of 
sectional curvature and diameter.
\end{abstract}

\section{Introduction} \label{sect1}

In a previous paper, we analyzed the limit of the $p$-form Laplacian
under a collapse with bounded sectional curvature and bounded diameter
to a smooth limit space \cite{Lott (2001)}.  In the present paper we
extend the analysis of \cite{Lott (2001)} to cover the case of a
singular limit space.

We give applications to upper and lower bounds on $\lambda_{p,j}$,
the $j$-th eigenvalue
of the $p$-form Laplacian counted with multiplicity, 
in terms of sectional curvature and diameter.
To do so, we use Gromov's precompactness theorem 
\cite[Chapter 5]{Gromov (1999)}.  Suppose that $M$ is a smooth connected
closed $n$-dimensional
manifold and $\{g_i\}_{i=1}^\infty$ is a sequence of Riemannian metrics
on $M$ of uniformly bounded sectional curvature and diameter.  Gromov's
theorem implies that a subsequence of the metric spaces 
$\{(M, g_i)\}_{i=1}^\infty$, which we relabel as $\{(M, g_i)\}_{i=1}^\infty$,
converges in the Gromov-Hausdorff topology to
a compact metric space $X$. If we can prove that (after passing to a
further subsequence) the $j$-th eigenvalue of the
$p$-form Laplacian on $(M, g_i)$ converges to the $j$-th eigenvalue of
an appropriate operator on $X$, and if we can effectively analyze the operator
on $X$, then we can obtain analytic results that are valid for all
Riemannian
metrics on $M$ with given bounds on sectional curvature and diameter.

In \cite{Lott (2001)} we showed that if
$X$ happens to be a smooth
Riemannian manifold $B$ then the relevant operator on $B$ is the Laplacian
associated to a flat degree-$1$ superconnection on $B$. In the general case,
the limit space $X$ may not be homeomorphic to a manifold.  However, Fukaya
showed that $X$ is the quotient
of a manifold by the action of a compact
Lie group \cite{Fukaya (1988)}. Namely,
a Riemannian metric $g$ on $M$ induces
a canonical Riemannian metric $\check{g}$ on the orthonormal frame bundle
$FM$ \cite[Section 1]{Fukaya (1988)}. 
The sectional curvature and diameter bounds on $\{g_i\}_{i=1}^\infty$
imply similar bounds on $\{\check{g}_i\}_{i=1}^\infty$. Then
$\{(FM, \check{g}_i)\}_{i=1}^\infty$ has a subsequence which converges to
a Riemannian
manifold $\check{X}$ in the $O(n)$-equivariant Gromov-Hausdorff topology,
and $X \: = \: \check{X}/O(n)$.

Hence one may hope to construct the
relevant operator on $X$ by working equivariantly on $\check{X}$.
This is in analogy to what Fukaya did in the case of the function
Laplacian \cite[Section 7]{Fukaya (1987)}. 

In this paper we carry out this program, using in part geometric results of
Cheeger, Fukaya and Gromov \cite{Cheeger-Fukaya-Gromov (1992)}. In effect
we construct superconnection Laplacians on the singular space $X$, and
prove convergence and compactness properties of the operators.
We now give some of the consequences.

\begin{definition} \label{def1}
If $M$ is a connected closed smooth manifold of positive dimension, 
let ${\cal M}(M,K)$ be the set of Riemannian metrics on $M$ with
$\diam(M) \: \le \: 1$ and $\parallel R^M \parallel_\infty \: \le \: K$.
Let ${\cal M}^1(M,K)$ be the subset with $\diam(M) \: = \: 1$.
Put
\begin{equation} \label{eq1.1} 
a_{p,j,K}(M) =  \inf_{{\cal M}(M,K)} \{ \lambda_{p,j}(M) \}, \: \:
\: \: \: \: \: \: \: \: \: \:  
a^1_{p,j,K}(M) =  \inf_{{\cal M}^1(M,K)} \{ \lambda_{p,j}(M) \},
\end{equation}
\begin{equation*}
A^1_{p,j,K}(M) =  \sup_{{\cal M}^1(M,K)} \{ \lambda_{p,j}(M) \}.
\end{equation*}
For $n \in \Z^+$, put
\begin{equation} \label{eq1.2} 
a_{n,p,j,K}^1 =  \inf_M \{ a_{p,j,K}^1(M) \},
\: \: \: \: \: \: \: \: \: \: 
A_{n,p,j,K}^1 = \sup_M \{ A_{p,j,K}^1(M) \},
\end{equation}
where $M$ ranges over connected closed smooth $n$-dimensional manifolds.
\end{definition}

Our first result says that there is a uniform lower bound on 
$\lambda_{p,j}$ in terms of sectional curvature and diameter. More
importantly, we also characterize when there is not a uniform upper bound.

\begin{theorem} \label{th5}
1. $\lim_{j \rightarrow \infty} a^1_{n,p,j,K} = \infty$.\\
2. If $p \: > \: 1$ and $A^1_{p,j,K}(M) \: = \: \infty$ then $M$ collapses
with bounded sectional curvature and bounded diameter to a limit space $X$
satisfying $1 \: \le \dim(X) \: \le p-1$.  Furthermore, there is a
(possibly singular) affine fiber bundle $M \rightarrow X$ whose
generic fiber $Z$ does not admit any nonzero affine-parallel
$k$-forms for $p - \dim(X) \: \le \: k \: \le \: p$.\\
3. If $p \in \{0,1\}$ then $A^1_{n,p,j,K} \: < \: \infty$.
\end{theorem}

Theorem \ref{th5}.1 also
follows from heat kernel bounds (\cite{Berard (1988)} and references
therein), with explicit estimates.  Our proof, which just
uses collapsing arguments, is given by way of illustration.
We remark that the Ricci-analog of Theorem \ref{th5}.1 is
false, as the Betti numbers can be arbitrarily large
\cite[Theorem 0.4]{Anderson (1992)}.
Theorem \ref{th5}.3 also follows from \cite{Cheng (1975)}, with
explicit estimates.

Theorem 1 of \cite{Lott (2001)} gives a partial
converse to Theorem \ref{th5}.2, 
in the sense that if $M$ is the total space of an affine fiber bundle 
\cite[Definition 1]{Lott (2001)} over
a smooth manifold $B$ with $1 \: \le \: \dim(B) \: \le \: p \: - \: 1$,
and the fiber $Z$ does not admit any
nonzero affine-parallel $k$-forms for $p - \dim(B) \: \le \: k \: \le \: p$,
then $A^1_{p,j,K}(M) \: = \: \infty$. 

Theorem \ref{th5}.2 implies 
the following characterization of manifolds $M$ which
do not have uniform upper bounds on the eigenvalues of the $2$-form
Laplacian.

\begin{theorem} \label{tcor1}
If $A^1_{2,j,K}(M) \: = \: \infty$ then \\
1. There is an almost flat manifold $Z$ which does not admit any nonzero
affine-parallel $1$-forms or $2$-forms, and an affine diffeomorphism $\phi \in
\Aff(Z)$ such that $M$ is diffeomorphic to the mapping torus of $\phi$, or \\
2. There are almost flat manifolds $Z$, $Z_1$, $Z_2$ which do not admit
any nonzero affine-parallel $1$-forms or $2$-forms, and surjective affine maps
$\phi_i : Z \rightarrow Z_i$ such that $M$ is homeomorphic to the 
double mapping cylinder $\cyl(\phi_1) \cup_Z \cyl(\phi_2)$.
\end{theorem}

The next result says that there are at most $\bb_1(M) \: + \: \dim(M)$ small
eigenvalues of the $1$-form Laplacian. It is an extension of
\cite[Corollary 1]{Lott (2001)}.

\begin{theorem} \label{tcor2} 
If $a_{1,j,K}(M) = 0$ then
$j \: \le \: \bb_1(M) \: + \: \dim(M)$.
More precisely, suppose that
$a_{1,j,K}(M) = 0$ and $j \: > \: \bb_1(M)$. Let $X$ be the limit space
coming from the collapsing argument. Then
\begin{equation} \label{ineq}
j \: \le \:
\bb_1(X) \: + \: \dim(M) \: - \: \dim(X) \: \le \: \bb_1(M) \: + \: \dim(M).
\end{equation} 
\end{theorem}
The first inequality in (\ref{ineq}) is sharp, for example, in the case
of the Berger collapse of $S^3$ to $S^2$ 
\cite[Example 1.2]{Colbois-Courtois (1990)}.

Next, we give 
a bound on the number of small eigenvalues of the $p$-form
Laplacian for a manifold which is Gromov-Hausdorff close to a codimension-$1$
space $X$.  It is an extension of
\cite[Corollary 2]{Lott (2001)} and 
\cite[Th\'eor\`eme 1.17]{Colbois-Courtois (2000)}.

\begin{theorem} \label{tcor4}
Let $X$ be a connected closed $(n-1)$-dimensional Riemannian
orbifold.  Then for any $K \ge 0$, there are
$\delta, c > 0$ with the following property :
Suppose that $M$ is a connected closed smooth $n$-dimensional 
Riemannian manifold with $\parallel R^M \parallel_\infty \: \le \: K$
and $d_{GH}(M,X) < \delta$. First, $M$ is the total space of
an orbifold circle bundle over $X$.  Let ${\cal O}$ be the
orientation bundle of $M \rightarrow X$, a flat real line bundle on
$X$ in the orbifold sense. Then
$\lambda_{p,j}(M, g) \: > \: c$ for 
$j \: = \: \bb_p(X) + \bb_{p-1}(X; {\cal O}) \: + \: 1$.
\end{theorem}

In \cite[Definition 3]{Lott (2001)} we defined the notion of a collapsing
sequence of metrics $\{g_i\}_{i=1}^\infty \subset {\cal M}(M, K)$ with a 
smooth limit space $B$. 
This means that there are an affine
fiber bundle $M \rightarrow B$
and an $\epsilon \: > \: 0$ such that
each $(M, g_i)$ is $\epsilon$-biLipschitz to a model metric on the 
total space of the affine fiber bundle. We remark that for a given
$\epsilon \: > \: 0$, results of Cheeger, Fukaya and Gromov imply that
if $(M, g) \in {\cal M}(M, K)$ is sufficiently close to $B$ in the 
Gromov-Hausdorff topology then $(M, g)$ is $\epsilon$-biLipschitz to
a model metric on some affine fiber bundle $M \rightarrow B$
\cite[Proposition 4.9]{Cheeger-Fukaya-Gromov (1992)}. Hence the content of
our assumption is that there is a single affine fiber bundle involved 
for all of the $g_i$'s.

There is an extension of the notion of a collapsing sequence to the case of
a singular limit space $X$. Namely,
a collapsing sequence consists of 
a sequence $\{g_i\}_{i=1}^\infty \subset {\cal M}(M, K)$ and
an $O(n)$-equivariant affine fiber bundle
$FM \rightarrow \check{X}$ such that
$\{\check{g}_i\}_{i=1}^\infty$ is a $O(n)$-equivariant
collapsing sequence on $FM$. Given
what is proven in this paper, the
results proved in \cite{Lott (2001)} concerning collapsing sequences, 
with smooth limit space $B$, extend to results concerning
collapsing sequences with limit
space $X$. We state one such result, which
is an extension of
\cite[Theorem 5]{Lott (2001)}. It says that there are
three mechanisms to make small positive eigenvalues of the
differential form Laplacian in a collapsing sequence.
Either the differential form Laplacian on the generic fiber of the map
$M \rightarrow X$ admits small positive
eigenvalues, or
the pushforward ``cohomology'' sheaf on $X$ fails to
be semisimple, or the Leray spectral sequence of the map
$M \rightarrow X$ does not
degenerate at the $E_2$-term.

\begin{theorem} \label{th6}
Let $\{(M, g_i)\}_{i=1}^\infty$ 
be a collapsing sequence with limit $X$.
Suppose that 
$\lim_{i \rightarrow \infty} \lambda_{p,j}(M, g_i) = 0$ for
some $j > \bb_p(M)$. 
Write the generic fiber $Z$ of the map $M \rightarrow X$
as the quotient of a
nilmanifold $\widehat{Z} \: = \: \widehat{\Gamma} \backslash N$ by a finite
group $F$.
Then \\
1. For some $q \in [0,p]$, 
$\bb_q(Z) \: < \: \dim \left(\Lambda^q({\frak n}^*)^F 
\right)$, or\\
2. For all 
$q \in [0,p]$, $\bb_q(Z) \: = \: \dim \left(\Lambda^q({\frak n}^*)^F 
\right)$, and for 
some $q \in [0, p]$, the ``cohomology'' sheaf
$\HH^q(Z; \R)$ on $X$ has a subsheaf which is not a direct summand, or\\
3. For all $q \in [0,p]$, 
$\bb_q(Z) \: = \: \dim \left(\Lambda^q({\frak n}^*)^F 
\right)$ and each subsheaf of the ``cohomology'' sheaf
$\HH^q(Z; \R)$ on $X$ is a direct summand, and
the Leray spectral sequence to compute $\HH^p(M; \R)$
does not degenerate at the $E_2$-term.
\end{theorem}

As one consequence of Theorem \ref{th6}, 
we obtain a characterization of when the $p$-form
Laplacian has small positive eigenvalues in a collapsing sequence over
a codimension-$1$ space. It is an extension of
\cite[Corollary 5]{Lott (2001)} and 
\cite[Th\'eor\`eme 1.17]{Colbois-Courtois (2000)}.

\begin{theorem} \label{tcor7}
Let $\{(M, g_i)\}_{i=1}^\infty$ 
be a collapsing sequence associated to a limit space $X$
with $\dim(X) = \dim(M) - 1$. Suppose that 
$\lim_{i \rightarrow \infty} \lambda_{p,j}(M, g_i) = 0$ for
some $j > \bb_p(M)$.  Let ${\cal O}$ be the
orientation bundle of the orbifold circle bundle
$M \rightarrow X$, a flat real line bundle on
the orbifold $X$ in the orbifold sense.
Let $\chi \in \HH^2(X; {\cal O})$ be the Euler class of the
orbifold circle bundle $M \rightarrow X$. Let ${\cal M}_\chi$ be multiplication
by $\chi$. Then
${\cal M}_\chi : 
\HH^{p-1}(X; {\cal O}) \rightarrow \HH^{p+1}(X; \R)$ is nonzero
or ${\cal M}_\chi : 
\HH^{p-2}(X; {\cal O}) \rightarrow \HH^{p}(X; \R)$ is nonzero.
\end{theorem}

The structure of the paper is as follows.
Section \ref{sect6} deals
with the construction of the Laplacian associated to a flat degree-$1$
superconnection on a singular space of the form $X = \check{X}/G$,
where $\check{X}$ is a smooth closed Riemannian manifold and
$G$ is a closed subgroup of $\Isom(\check{X})$. 
In Section \ref{sect7} we prove
Theorems \ref{th5} and
\ref{tcor1}. Section \ref{sect8} 
uses the compactness results to prove
Theorems \ref{tcor2}-\ref{tcor7}.
More detailed descriptions appear at the
beginnings of the sections

As the present paper is a sequel to \cite{Lott (2001)}, we
sometimes make reference to the relevant sections of \cite{Lott (2001)}.
As for notation in this paper, 
if $G$ is a group which acts on a set $X$, we let $X^G$ denote the
set of fixed-points. If $B$ is a smooth manifold and $E$ is a smooth
vector bundle on $B$, we let $\Omega(B; E)$ denote the smooth
$E$-valued differential forms on $B$. If ${\frak n}$ is a nilpotent
Lie algebra on which a finite group $F$ acts by automorphisms
then ${\frak n}^*$ denotes the dual space, $\Lambda^*({\frak n}^*)$ denotes
the exterior algebra of the dual space and $\Lambda^*({\frak n}^*)^F$ denotes
the $F$-invariant subspace of the exterior algebra.

\section{Basic Laplacian} \label{sect6}

In this section we construct differential form Laplacians on a certain class
of singular spaces, namely those of the form $X = \check{X}/G$ where
$\check{X}$ is a smooth closed Riemannian manifold and $G$ is a closed
subgroup of $\Isom (\check{X})$.
Let $\check{E}$ be a $G$-equivariant vector bundle on $\check{X}$. We 
consider the space of basic $\check{E}$-valued
forms $\Omega_{basic}(\check{X}; \check{E})$.
If $\check{A}^\prime$ is a basic flat degree-$1$ superconnection on 
$\check{E}$ then we construct the corresponding Laplacian $\triangle^E$ as an
operator on $\Omega_{basic}(\check{X}; \check{E})$.
Although it is not strictly necessary for this paper,
we also give a more intrinsic formulation of $\triangle^E$ as an operator
on a space $\Omega(X; E)$ of forms on the quotient space $X$. We then
describe a spectral sequence to compute the cohomology of $\check{A}^\prime$.

We remark that the spaces $X = \check{X}/G$ can be quite singular.  
For example,
if $G$ is finite then one finds the well-known orbifold singularities, 
whereas if
$G$ has positive dimension then the singularities of $X$ can be much worse.
In view of the well-known difficulties in doing analysis on singular spaces, 
one may wonder how we can construct reasonable operators 
on such spaces.
The point is that there are special features of the present situation.
For example, there is an induced measure on $X$, the pushforward measure
from $\check{X}$, which has the effect of
mollifying the singularities.

Let $\check{X}$ be a smooth connected 
closed Riemannian manifold on which a compact
Lie group $G$ acts isometrically on the right. 
Let ${\frak g}$ denote the Lie algebra of $G$.

The $G$-action partitions $\check{X}$ into smooth submanifolds 
$\check{X}^{[H]}$,
where $[H]$ runs over the conjugacy classes of closed subgroups of $G$, and 
where $\check{x} \in \check{X}^{[H]}$ 
if the isotropy subgroup $G_{\check{x}}$
is in the conjugacy class $[H]$. The (connected components of the)
$\check{X}^{[H]}$'s
induce a stratification of the
quotient space $X = \check{X}/G$ as
$X = \bigcup_{[H]} X^{[H]}$, where $X^{[H]} = \check{X}^{[H]}/G$. The
orbits of $G$ on $\check{X}^{[H]}$ are all diffeomorphic to $H \backslash G$.
In fact, $\check{X}^{[H]}$ is a fiber bundle over ${X}^{[H]}$ with fiber
$H \backslash G$. To describe the structure group of this fiber bundle,
note that the ``internal symmetry group'' 
$\left( \Diff(H \backslash G) \right)^G$ of a fiber
$H \backslash G$ is isomorphic to $H \backslash N_H(G)$,
where $N_H(G)$ is the normalizer of
$H$ in $G$. (An element $n \in N_H(G)$ sends $Hg$ to $Hng$.)  Then the
structure group of the fiber bundle 
$\check{X}^{[H]} \rightarrow {X}^{[H]}$
is contained in $H \backslash N_H(G)$.
To see this more explicitly, put $\check{X}^{(H)} \: = \: 
\{\check{x} \in \check{X} \: : \: G_{\check{x}} \: = \: H\}$.
Then $H \backslash N_H(G)$ acts freely on 
$\check{X}^{(H)}$.
In fact, $\check{X}^{(H)}$ is a principal 
$H \backslash N_H(G)$-bundle. Then
$\check{X}^{[H]} \: = \: \check{X}^{(H)} \times_{H \backslash N_H(G)} 
(H \backslash G)$ and 
${X}^{[H]} \: = \: \check{X}^{(H)}/(H \backslash N_H(G))$.

There is another stratification of $X$, introduced in 
\cite{Davis (1978)}, which is more convenient for our purposes. It 
keeps track of both the connected components of $\check{X}^{[H]}$
and their normal bundles in $\check{X}$.
We briefly review the setup of \cite{Davis (1978)}. 
Consider pairs $(H, V)$ where $H$ is a closed subgroup of
$G$ and $V$ is a real representation space of $H$ with no trivial
subrepresentations, i.e. $V^H = \{0\}$. There is a natural $G$-action on
such pairs and the equivalence classes are called the normal orbit types 
$[H,V]$.
Given a point $\check{x} \in \check{X}$, the differentiable
slice theorem says that there is a real representation space $W_{\check{x}}$
of $G_{\check{x}}$ so that a neighborhood of $\check{x} \cdot G$ is
$G$-diffeomorphic to $W_{\check{x}} \times_{G_{\check{x}}} G$. 
Put 
$V_{\check{x}} = W_{\check{x}} / W_{\check{x}}^{G_{\check{x}}}$.
Then $[G_{\check{x}}, V_{\check{x}}]$ is called the normal orbit type
of $\check{x}$. Given a normal orbit type $\alpha$, put
\begin{equation} \label{eq6.1}
\check{X}_\alpha \: = \: \left\{ \check{x} \in \check{X} :
[G_{\check{x}}, V_{\check{x}}] = \alpha \right\}.
\end{equation}
Then $\check{X}_\alpha$ is a smooth $G$-submanifold of $\check{X}$.
Let ${\cal N}$ be the set of normal orbit types $\alpha$ such that 
$\check{X}_\alpha$ is nonempty. For $\alpha \in {\cal N}$, put
$X_\alpha = \check{X}_\alpha / G$, a smooth Riemannian manifold. 
Then $X$ is stratified by 
$\{X_\alpha\}_{\alpha \in {\cal N}}$.
Once again, $\check{X}_\alpha$ is a fiber bundle over
$X_\alpha$ with fiber $H \backslash G$.

As $G$ acts isometrically, we may assume that for each normal orbit type
$\alpha = [H, V]$, $V$ is given an $H$-invariant
inner product. Let $SV$ denote the corresponding sphere in $V$.

Given a normal orbit type $\alpha$, fix a representative $(H, V)$. Consider
the diagonal inclusion of $H$ in $O(V) \times G$. Let
$N_H (O(V) \times G)$ be the normalizer of $H$ in $O(V) \times G$.
Let $\nu_\alpha$ be the normal bundle of $\check{X}_\alpha$ in
$\check{X}$. The restriction of $\nu_\alpha$ to a fiber
$H \backslash G$ of $\check{X}_\alpha \rightarrow X_\alpha$ is
isomorphic to 
the Euclidean vector bundle 
$(V \times_H G) \rightarrow (H \backslash G)$.
The ``internal symmetry group'' of this vector bundle,
i.e. the group of vector bundle automorphisms
which commute with the $G$-action, is 
$S_\alpha \equiv H \backslash N_H (O(V) \times G)$.
Correspondingly, there is
a certain principal $S_\alpha$-bundle $P_\alpha$ such that 
\begin{equation} \label{eq6.2}
\nu_\alpha = P_\alpha \times_{S_\alpha} (V \times_H G).
\end{equation}
Using the normal exponential map, 
there is a neighborhood of $\check{X}_\alpha$ in
$\check{X}$ which is $G$-diffeomorphic to $\nu_\alpha$.
In addition,
\begin{equation}
\check{X}_\alpha = P_\alpha \times_{S_\alpha} (H \backslash G),
\end{equation}
\begin{equation}
X_\alpha = P_\alpha / S_\alpha
\end{equation}
and there is a 
neighborhood $N_\alpha$ of $X_\alpha$ in $X$ whose homeomorphism type is
\begin{equation} \label{eq6.3}
N_\alpha \: \cong \: P_\alpha \times_{S_\alpha} (V /H).
\end{equation}
(Note that $N_\alpha$ and $V/H$ are generally not manifolds.)

There is a partial ordering on normal orbit types given by saying that
$[H,V] \le [H^\prime,V^\prime]$ if the $G$-manifold $V \times_H G$ contains
a $G$-orbit of type $[H^\prime, V^\prime]$. This induces a partial ordering
on the strata of $X$, with $\alpha \le \alpha^\prime$ if and only if
$X_{\alpha} \subset \overline{X_\alpha^\prime}$. 

Let $\check{E} = \oplus_{j=0}^m \check{E}^j$ be a 
$\Z$-graded real vector bundle on $\check{X}$.
We assume that
the action of $G$ on $\check{X}$ is covered by an action on $\check{E}$ 
which preserves the $\Z$-grading. Let ${\frak g}$ be the Lie algebra of $G$.
We say that $\check{E}$ is $G$-basic if it is equipped with a 
$G$-equivariant linear map
${\frak I} \: : \: {\frak g} \rightarrow 
C^\infty \left( \check{X}; \Hom(\check{E}^*, \check{E}^{*-1})  \right)$
which satisfies ${\frak I}({\frak x})^2 \: = \: 0$ for all
${\frak x} \in {\frak g}$.
Given ${\frak x} \in {\frak g}$, we write
${\frak I}_{\frak x}$ for ${\frak I}({\frak x})$.

Given 
${\frak x} \in {\frak g}$, let ${\frak X}$ be the corresponding vector field
on $\check{X}$ and let $i_{\frak X}$ be interior multiplication by
${\frak X}$ on $\Omega \left( \check{X}; \check{E} \right)$.

\begin{definition} \label{def4}
Put 
\begin{equation} \label{eq6.4}
\Omega_{G} \left( \check{X}; \check{E} \right) \: = \:
\{ \check{\omega} \in  \Omega \left( \check{X}; \check{E} \right) \: : \:
g \cdot \check{\omega} = \check{\omega} \text{ for all $g \in G$}\}
\end{equation}
and
\begin{equation} \label{eq6.5}
\Omega_{basic} \left( \check{X}; \check{E} \right) \: = \:
\{ \check{\omega} \in  \Omega_G \left( \check{X}; \check{E} \right) \: : \:
(i_{\frak X} \: + \: {\frak I}_{\frak x}) 
\check{\omega} \: = \: 0 \text{ for all ${\frak x} \in {\frak g}$}\}.
\end{equation}
\end{definition}

We give $\Omega_{basic} \left( \check{X}; \check{E} \right)$ the total
$\Z$-grading coming from the $\Z$-gradings on 
$\Omega^* ( \check{X} )$ and $\check{E}$.
Similarly, we let $\Omega_{basic} \left( \check{X}; \End(\check{E}) \right)$ be
the $G$-invariant elements $\check{\eta}$ of 
$\Omega \left( \check{X}; \End(\check{E}) \right)$ which satisfy
$i_{\frak X} \: \check{\eta} \: + \: {\frak I}_{\frak x} \:
\check{\eta} \: + \: (-1)^{|\check{\eta}|} \: \check{\eta} \:
{\frak I}_{\frak x} \: = \: 0$.

Let $h^{\check{E}}$ be a $G$-invariant graded 
Euclidean inner product on $\check{E}$.  We obtain 
$L^2$-inner products on $\Omega_{G} \left( \check{X}; 
\check{E} \right)$ and
$\Omega_{basic} \left( \check{X}; \check{E} \right)$.
Let $\Omega_{G, L^2} \left( \check{X}; \check{E} \right)$
and $\Omega_{basic, L^2} \left( \check{X}; \check{E} \right)$ be the 
$L^2$-completions of 
$\Omega_{G} \left( \check{X}; \check{E} \right)$ and
$\Omega_{basic} \left( \check{X}; \check{E} \right)$, respectively.
Let $P^{hor} : \Omega_{G, L^2} 
\left( \check{X}; \check{E} \right) \rightarrow
\Omega_{basic, L^2} \left( \check{X}; \check{E} \right)$ be 
orthogonal projection.

\begin{definition} \label{def5}
Let $\check{E}$ be a real 
$G$-basic vector bundle on $\check{X}$.
A basic connection on $\check{E}$ is
a connection $\nabla^{\check{E}}$ on $\check{E}$ 
which is $G$-invariant
and satisfies the property that
for all ${\frak x} \in {\frak g}$, 
$\nabla^{\check{E}} (i_{\frak X} \: + \: {\frak I}_{\frak x}) \: + \:
(i_{\frak X} \: + \: {\frak I}_{\frak x}) \nabla^{\check{E}}$ 
equals Lie differentiation ${\cal L}_{\frak X}$ with respect to ${\frak X}$ on
$\Omega \left( \check{X}; \check{E} \right)$.
\end{definition}

If $\nabla^{\check{E}}$ is a basic connection then
$\nabla^{\check{E}} i_{\frak X} \: + \:
i_{\frak X} \nabla^{\check{E}} \: = \:
{\cal L}_{\frak X}$ and 
$\nabla^{\check{E}} {\frak I}_{\frak x} \: + \:
{\frak I}_{\frak x} \nabla^{\check{E}} \: = \: 0$, i.e.
$\nabla^{\check{E}}$ is basic in the usual sense and 
${\frak I}_{\frak x}$ is covariantly-constant with respect to
$\nabla^{\check{E}}$.
A basic flat connection is a connection which is both basic and flat.

We wish to describe the set of basic flat connections on $\check{X}$ in terms
of a representation variety.  To do so, we need the 
correct analog of the fundamental
group of $\check{X}/G$. Let $\widehat{X}$ be the universal cover of 
$\check{X}$, with
projection $q : \widehat{X} \rightarrow \check{X}$. Put 
\begin{equation} \label{eq6.6}
\widehat{G} = \{ (\phi, g) \in \Diff(\widehat{X}) \times G \: : \: 
q \circ \phi \: = \: g \cdot q \}.
\end{equation}
There is an exact sequence of groups
\begin{equation} \label{eq6.7}
1 \longrightarrow \pi_1(\check{X})
\longrightarrow \widehat{G}
\longrightarrow G \longrightarrow 1.
\end{equation}
The corresponding homotopy exact sequence of spaces gives
\begin{equation} \label{eq6.8}
\ldots \rightarrow
1 \rightarrow \pi_1(\widehat{G}) \rightarrow \pi_1(G) \rightarrow
\pi_1(\check{X}) 
\rightarrow \pi_0(\widehat{G}) \rightarrow \pi_0(G)
\rightarrow 1.
\end{equation}

The next proposition is implicit in
\cite[Section 4a]{Bismut-Lott (1995)}.

\begin{proposition} \label{prop6}
Suppose that for all ${\frak x} \in {\frak g}$, 
${\frak I}_{\frak x} \: = \: 0$.
Then there is a bijection between 
$\Hom \left( \pi_0(\widehat{G}), \GL(N, \R) \right)/\GL(N, \R)$ and
the basic flat connections on rank-$N$ $G$-vector bundles over
$\check{X}$, up to $G$-equivariant
gauge equivalence.
\end{proposition}
\begin{pf}
Given $\rho \in \Hom \left( \pi_0(\widehat{G}), \GL(N, \R) \right)$,
let $\overline{\rho}$ be the restriction of $\rho$ to 
$\pi_1(\check{X})$. 
Let $\rho_0 : \widehat{G} \rightarrow \pi_0( \widehat{G}) 
\stackrel{\rho}{\rightarrow} \GL(N, \R)$
be the composite homomorphism. Put $\check{E} = \widehat{X} 
\times_{\overline{\rho}} \R^N$. 
Then $\check{E}$ is a flat vector bundle on $\check{X}$.
If $\widehat{x} \in \widehat{X}$,
$v \in \R^N$ and $\widehat{g} \in \widehat{G}$, put
\begin{equation} \label{eq6.9}
(\widehat{x}, v) \cdot \widehat{g} = (\widehat{x} \cdot \widehat{g},
\rho_0(\widehat{g}^{-1}) v).
\end{equation}
Then  this action of $\widehat{G}$ on $\widehat{X} \times \R^N$ extends
that of the normal subgroup $\pi_1(\check{X})$ of $\widehat{G}$. 
Hence the quotient group $G$ acts on $\check{E}$. 
As the representation $\rho_0$ factorizes through $\rho$, we see that the flat
connection on $\check{E}$ is basic. Conjugate representations $\rho$ give
gauge-equivalent basic flat connections.

Conversely,
let $\nabla^{\check{E}}$ be a basic flat connection on a rank-$N$ $G$-vector
bundle $\check{E}$. Putting $\widehat{E} = q^* \check{E}$,
there is a trivialization $\widehat{E} \cong \widehat{X} \times \R^N$.
The action of $G$ on $\check{E}$ lifts to an action
of $\widehat{G}$ on $\widehat{E}$. In this way we get a
homomorphism $\rho_0 : \widehat{G} \rightarrow \GL(N, \R)$. As
$\nabla^{\check{E}}$ is basic, it follows that $\rho_0$ factors through
a representation $\rho : \pi_0(\widehat{G}) \rightarrow \GL(N, \R)$.
Taking into account the ambiguity in the trivialization of $\widehat{E}$,
we obtain a well-defined
conjugacy class of $\rho$.
\end{pf} 

If ${\frak I}_{\frak x}$ is not identically zero 
then we can describe the basic flat
connections as the subset of the connections in Proposition \ref{prop6} 
with respect to which
${\frak I}_{\frak x}$ is covariantly-constant for all
${\frak x} \in {\frak g}$.

For background information on superconnections, we refer to
\cite[Chapter 1.4]{Berline-Getzler-Vergne (1992)}.
\begin{definition} \label{def6}
A basic superconnection on $\check{E}$ is
a superconnection $\check{A}^\prime$ on $\check{E}$ which is $G$-invariant
and satisfies the property that
for all ${\frak x} \in {\frak g}$, 
$\check{A}^\prime (i_{\frak X} + {\frak I}_{\frak x}) \:  + \:
(i_{\frak X} + {\frak I}_{\frak x}) \check{A}^\prime$ 
equals Lie differentiation with respect to ${\frak X}$ on
$\Omega \left( \check{X}; \check{E} \right)$.
\end{definition}

A basic superconnection on $\check{E}$
restricts to a superconnection on 
$\Omega_{basic} \left( \check{X}; \check{E} \right)$.
Let $\check{A}^\prime$ be a basic superconnection on
$\check{E}$ which is ``flat degree-$1$''
in the sense of \cite[Section II(a)]{Bismut-Lott (1995)} and
\cite[Section 5]{Lott (2001)}.
Let $\left( \check{A}^\prime \right)^*$ be its formal adjoint on
$\Omega_{G} \left( \check{X}; \check{E} \right)$, which is the restriction
of the formal adjoint on $\Omega \left( \check{X}; \check{E} \right)$.

\begin{lemma} \label{lemma5}
The formal adjoint of $\check{A}^\prime$ on 
$\Omega_{basic} \left( \check{X};
\check{E} \right)$ is 
\begin{equation} \label{eq6.11}
\left( \check{A}^\prime \right)^*_{basic} \: = \: P^{hor} 
\left( \check{A}^\prime \right)^*.
\end{equation}
\end{lemma}
\begin{pf}
Given $\check{\omega}, \check{\omega}^\prime 
\in \Omega_{basic} \left( \check{X};
\check{E} \right)$, we have
\begin{equation} \label{eq6.12}
\langle \check{\omega}, \check{A}^\prime \check{\omega}^\prime \rangle \: = \:
\langle \left( \check{A}^\prime \right)^* \check{\omega}, 
\check{\omega}^\prime \rangle \: = \:
\langle \left( \check{A}^\prime \right)^* \check{\omega}, P^{hor} 
\check{\omega}^\prime \rangle \: = \:
\langle  P^{hor} \left( \check{A}^\prime \right)^* \check{\omega}, 
\check{\omega}^\prime \rangle.
\end{equation}
As $P^{hor} \left( \check{A}^\prime \right)^* \check{\omega} \in 
\Omega_{basic, L^2} \left( \check{X}; \check{E} \right)$, the lemma
follows.
\end{pf}

Put $\Omega_{basic, max}(\check{X}; \check{E}) \: = \: 
\left\{\check{\omega} \in \Omega_{basic,L^2}(\check{X}; \check{E}) \: : \:
\check{A}^\prime \check{\omega} \in \Omega_{basic,L^2}(\check{X}; \check{E})
\right\} \subset \Omega_{basic,L^2}(\check{X}; \check{E})$,
where $\check{A}^\prime \check{\omega}$ 
is originally defined as a distribution.

\begin{lemma}
The operator $\check{A}^\prime$ is closed on 
$\Omega_{basic, max}(\check{X}; \check{E})$.
\end{lemma}
\begin{pf}
Suppose that $\{\check{\omega}_i\}_{i=1}^\infty$ is a sequence in
$\Omega_{basic, max}(\check{X}; \check{E})$ such that
$\lim_{i \rightarrow \infty} \check{\omega}_i \: = \: \check{\omega}$ in
$\Omega_{basic, L^2}(\check{X}; \check{E})$ and
$\lim_{i \rightarrow \infty} \check{A}^\prime \check{\omega}_i \: = \: 
\check{\eta}$ in
$\Omega_{basic, L^2}(\check{X}; \check{E})$.
Given $\check{\phi} \in \Omega \left( \check{X}; \check{E} \right)$, we have
\begin{equation}
\langle \check{\phi}, \check{\eta} \rangle \: = \: \lim_{i \rightarrow \infty}
\langle \check{\phi}, \check{A}^\prime \check{\omega}_i \rangle \: = \: 
\lim_{i \rightarrow \infty} \langle (\check{A}^\prime)^* \check{\phi}, 
\check{\omega}_i \rangle \: = \: \langle (\check{A}^\prime)^* \check{\phi}, 
\check{\omega} \rangle.
\end{equation}
It follows that $\check{\omega} 
\in \Omega_{basic, max}(\check{X}; \check{E})$ with
$\check{A}^\prime \: \check{\omega} \: = \: \check{\eta}$.
This proves the lemma.
\end{pf}

Let $\left( \check{A}^\prime \right)^*_{basic}$ be the adjoint operator
to $\check{A}^\prime$, the latter being defined on 
$\Omega_{basic, max}(\check{X}; \check{E})$.
As $\Image \left( \check{A}^\prime \right) \subset
\Ker \left( \check{A}^\prime \right)$, a formal result of functional
analysis \cite[Lemma 4.3]{Hilsum (1985)}
implies that there is an orthogonal decomposition
\begin{equation} \label{ortho}
\Omega_{basic, L^2}(\check{X}; \check{E}) \: = \:
\left( \Ker(\check{A}^\prime) \cap \Ker((\check{A}^\prime)^*_{basic}) \right)
\oplus 
\frac{\Omega_{basic, L^2}(\check{X}; \check{E})}{\Ker(\check{A}^\prime)} 
\oplus 
\frac{\Omega_{basic, L^2}(\check{X}; 
\check{E})}{\Ker((\check{A}^\prime)^*_{basic})}.
\end{equation}
Furthermore, 
$\check{A}^\prime \: + \: (\check{A}^\prime)^*_{basic}$ is self-adjoint on
$\Omega_{basic, L^2}(\check{X}; \check{E})$.

\begin{definition} \label{def7}
The basic Laplacian is
$\triangle^{\check{E}}_{basic} \: = \: 
\left( \check{A}^\prime \: + \: (\check{A}^\prime)^*_{basic}) \right)^2$.
\end{definition}

In terms of (\ref{ortho}), the domain of $\triangle^{\check{E}}_{basic}$ is
\begin{align}
\left( \Ker(\check{A}^\prime) \cap \Ker((\check{A}^\prime)^*_{basic}) \right)
& \oplus \frac{\{ \check{\omega}  \in 
\Dom(\check{A}^\prime) \: : \: \check{A}^\prime \check{\omega} \in 
\Dom((\check{A}^\prime)^*_{basic}) \}}{\Ker(\check{A}^\prime)} \\
& \oplus \frac{\{ \check{\omega}  \in 
\Dom((\check{A}^\prime)^*_{basic}) 
\: : \: (\check{A}^\prime)^*_{basic} \check{\omega} \in 
\Dom(\check{A}^\prime) \}}{\Ker((\check{A}^\prime)^*_{basic})}, \notag
\end{align}
on which $\triangle^{\check{E}}_{basic}$
acts by $0 \: + \: (\check{A}^\prime)^*_{basic} \check{A}^\prime \: + \:
\check{A}^\prime (\check{A}^\prime)^*_{basic}$. 
The quadratic form $Q$ associated
to $\triangle^{\check{E}}_{basic}$ has domain
\begin{equation}
\Dom(Q) \: = \: \{ \check{\omega} \in \Omega_{basic, L^2}(\check{X}; \check{E})
\: : \: \check{A}^\prime \check{\omega} \in
\Omega_{basic, L^2}(\check{X}; \check{E}) \text{ and }
(\check{A}^\prime)^*_{basic} \check{\omega} \in
\Omega_{basic, L^2}(\check{X}; \check{E}) \},
\end{equation}
on which it is defined by
\begin{equation}
Q(\check{\omega}) \: = \:
\langle \check{A}^\prime \check{\omega}, 
\check{A}^\prime \check{\omega} \rangle \: + \:
\langle (\check{A}^\prime)^*_{basic} \check{\omega}, 
(\check{A}^\prime)^*_{basic} \check{\omega} \rangle.
\end{equation}

We can also describe $\triangle^{\check{E}}_{basic}$ when restricted to
$\Image(\check{A}^\prime)^\perp$ in terms of a
Friedrichs extension.   Namely, using Lemma \ref{lemma5}, we have that
$(\check{A}^\prime)^*_{basic} \check{A}^\prime$ maps 
$\Omega_{basic}(\check{X}; 
\check{E})$ to $\Omega_{basic, L^2}(\check{X}; 
\check{E})$. Then the Friedrichs extension
\cite[Theorem X.23]{Reed-Simon (1978)} of the
positive symmetric operator 
$(\check{A}^\prime)^*_{basic} \check{A}^\prime$, acting on  
$\Omega_{basic}(\check{X}; 
\check{E})$, is well-defined, coincides with
$\triangle^{\check{E}}_{basic}$
on 
\begin{equation}
\Image(\check{A}^\prime)^\perp \: = \: 
\left( \Ker(\check{A}^\prime) \cap \Ker((\check{A}^\prime)^*_{basic}) \right)
\oplus 
\frac{\Omega_{basic, L^2}(\check{X}; \check{E})}{\Ker(\check{A}^\prime)}
\end{equation}
and vanishes on
$\frac{\Omega_{basic, L^2}(\check{X}; 
\check{E})}{Ker((\check{A}^\prime)^*_{basic})}$.

Using the orthogonal decomposition (\ref{ortho}), we have that
$\Ker(\triangle^{\check{E}}_{basic}) \: = \: \Ker(\check{A}^\prime)/
\overline{\Image(\check{A}^\prime)}$. If $\check{A}^\prime$ has a closed
image then $\Ker(\triangle^{\check{E}}_{basic})$ equals the usual
cohomology $\Ker(\check{A}^\prime)/
{\Image(\check{A}^\prime)}$. We do not show that 
$\check{A}^\prime$ has a closed
image in general, but we will see in Section \ref{sect8} that it has
a closed image in a case that arises in collapsing.

Let us look in more detail at the basic forms when $\check{X} \: = \:
H \backslash G$. If $\check{E}_{He}$ denotes the fiber of $\check{E}$ at
$He \in H \backslash G$ then
we can write the homogeneous vector bundle $\check{E}$ as
$\check{E} \: = \: \check{E}_{He} \times_H G$. The restriction of
${\frak I}$ to $\check{E}_{He}$ becomes an $H$-equivariant map
${\frak I} \: : \: {\frak g} \rightarrow 
\Hom \left( \check{E}_{He}^*, \check{E}_{He}^{*-1} \right)$.
Let ${\frak h}$ be the Lie algebra of $H$. 
Put
$\check{K} \: = \: \bigcap_{x \in {\frak h}} 
\Ker \left( {\frak I}_{\frak x} \right)$.
Let $\{ {\frak x}_j \}$ be a basis of
${\frak g}/{\frak h}$ and let
$\{ {\frak x}_j^* \}$ be the dual basis of
$({\frak g}/{\frak h})^*$. 
If $\check{\omega}$ is a basic form
then as $\left( (i_{\frak X} \: + \: {\frak I}_x) \check{\omega} \right) (He)
\: = \: 0 $ for all ${\frak x} \in {\frak g}$, it follows that
$\check{\omega} (He) \: = \: \prod_j (1 \: - \: e({\frak x}_j^*) \:
{\frak I}_{{\frak x}_j}) \: \check{\rho}$ for some
$\check{\rho} \in \check{K}^H$. Here $e({\frak x}_j^*)$ denotes exterior
multiplication by ${\frak x}_j^*$ and 
${\frak I}_{{\frak x}_j} \check{\rho}$
is well-defined since $\check{\rho} \in \check{K}$. The value of
$\check{\omega}$ on the rest of $H \backslash G$ is determined by
$G$-invariance.  In this way, we see that
$\Omega_{basic}(H \backslash G; \check{E}) \cong \check{K}^H$.

To extend this to general $\check{X}$,
suppose that 
$\check{x} \in \check{X}_\alpha$ has isotropy group
$H$, with Lie algebra ${\frak h}$.
Put
\begin{equation}
\check{K}_{\check{x}} \: = \:
\bigcap_{{\frak x} \in {\frak h}}
\Ker \left( {\frak I}_{\frak x}(\check{x}) \right).
\end{equation}
{\it A priori}
$\{\check{K}_{\check{x}}\}_{x \in \check{X}_\alpha}$ 
may not form a vector bundle
on $\check{X}_\alpha$, due to possible jumps in the dimension. 
Hereafter we make the assumption that for each closed subgroup $H$ of
$G$, $\left\{ \bigcap_{{\frak x} \in {\frak h}}
\Ker \left( {\frak I}_{\frak x}(\check{x}) \right) \right\}_{\check{x} \in
\check{X}^H}$ forms a vector bundle
on $\check{X}^H$; this is the case that will arise in collapsing.
It then follows from $G$-equivariance that
$\{\check{K}_{\check{x}}\}_{x \in \check{X}_\alpha}$ 
forms a vector bundle ${\check K}_\alpha$ on $\check{X}_\alpha$.
By $G$-equivariance, $\check{K}_\alpha$
gives rise to a $\Z$-graded real vector bundle $E_\alpha$ on $X_\alpha$
such that $C^\infty \left(X_\alpha; E_\alpha \right) \cong 
C^\infty \left( \check{X}_\alpha; \check{K}_\alpha \right)^G$.
We will sometimes write $E$ as shorthand
for $\{E_\alpha\}_{\alpha \in {\cal N}}$.
We remark that $\dim(E_\alpha)$ may not be constant in $\alpha$.
However, we will see that if $\alpha < \alpha^\prime$ then
$\dim(E_\alpha) \le \dim(E_{\alpha^\prime})$.

Given $\check{\omega} \in \Omega_{basic} \left( \check{X}; \check{E} \right)$,
we obtain a collection of forms 
$\left\{\check{\omega}_\alpha \in \Omega_{basic} \left(
\check{X}_\alpha; \check{E}_\alpha \right)\right\}_{\alpha \in {\cal N}}$ 
by pullback to the $\check{X}_\alpha$'s.
 Let $\check{\rho}_\alpha$ be the
horizontal component of $\check{\omega}_\alpha$ with respect to the
fiber bundle $\check{X}_\alpha \rightarrow X_\alpha$. That is, if
$\check{x} \in \check{X}_\alpha$ has isotropy group $H$,
let $\{ {\frak x}_j \}$ be a basis of
${\frak g}/{\frak h}$ and let $\{ {\frak x}_j^* \}$ be the dual basis of
$({\frak g}/{\frak h})^*$. Then 
\begin{equation} \label{to1}
\check{\rho}_\alpha(\check{x}) \: = \:
\prod_j (1 \: - \: e({\frak x}_j^*) \:
i_{{\frak X}_j}) \: \check{\omega}_\alpha
(\check{x}),
\end{equation}
where $e({\frak x}_j^*)$ denotes exterior multiplication by the
$1$-form in $T^*_{\check{x}} \check{X}_\alpha$
represented by ${\frak x}_j^*$ which is vertical with respect to
the Riemannian fiber bundle
$H \backslash G \rightarrow
\check{X}_\alpha \rightarrow X_\alpha$.
By construction, $\check{\rho}_\alpha$ is a $G$-invariant
element of $\Omega \left( \check{X}_\alpha; {\check K}_\alpha \right)$ and
satisfies $i_{\frak X} \check{\rho}_\alpha \: = \: 0$ for all
${\frak x} \in {\frak g}$.

We obtain a collection of forms
$\left\{\omega_\alpha \in \Omega \left(
X_\alpha; E_\alpha \right) \right\}_{\alpha \in {\cal N}}$ 
by subsequent pushforward of the $\check{\rho}_\alpha$'s to the $X_\alpha$'s. 
We define elements of  $\Omega \left( X; E \right)$ to be collections which
so arise. Note that given $\check{\rho}_\alpha(\check{x})$, we can recover
$\check{\omega}_\alpha(\check{x})$ by
\begin{equation} \label{rhs}
\check{\omega}_\alpha(\check{x}) \: = \: \prod_j (1 \: - \:
e({\frak x}_j^*) \: {\frak I}_{{\frak x}_j}) \: \check{\rho}_\alpha(\check{x}).
\end{equation}
In this way, there is an isomorphism between
$\Omega_{basic} (\check{X}; \check{E})$ and $\Omega (X; E)$.
If $\check{E}$ is the trivial $\R$-bundle on $\check{X}$ then we
denote $\Omega (X; E)$ by $\Omega^* (X)$.

Define $v_\alpha : X_\alpha \rightarrow \R^+$ by saying that for
$x \in X_\alpha$, $v_\alpha(x)$ is the volume in 
$\check{X}_\alpha$
of the $G$-orbit corresponding to $x$.
Then there are inner products
$\left\{ h^{E_\alpha} \right\}_{\alpha \in {\cal N}}$ on the $E_\alpha$'s
so that the isomorphism
from $\Omega_{basic} (\check{X}_\alpha; \check{E}_\alpha)$ 
to $\Omega (X_\alpha; E_\alpha)$ becomes an isometry, where
the inner product on $\Omega (X_\alpha; E_\alpha)$ is weighted by
$v_\alpha$.
If $\beta$ is the normal orbit type of the principal part of $\check{X}$ then
define a new inner product
$h^E$ on $E_\beta$ by $h^E = v_\beta \cdot h^{E_\beta}$.
Let $\Omega_{L^2} (X; E)$ be the
$L^2$-closure of $\Omega \left( X; E \right)$, 
using the Riemannian metric on $X_\beta$ and $h^E$.
There is an isometric isomorphism between
$\Omega_{basic, L^2} \left( \check{X}; \check{E} \right)$ and
$\Omega_{L^2}(X; E)$, so we can think of $\triangle^{\check{E}}_{basic}$ as
an operator $\triangle^E$ which is densely-defined on $\Omega_{L^2} (X; E)$.

In the special case when $\alpha$ is a normal orbit type with 
$V \: = \: \R$, we will need an additional vector bundle on
$X_\alpha$. (Note that this happens when $\check{X_\alpha}$ has codimension
one in $\check{X}$.) If $\alpha \: = \: [H, \R]$ then the representation
$H \rightarrow O(\R) ( \: = \: \Z_2)$ is necessarily onto, with kernel
conjugate to $H_\beta$. Let
$\nu_\alpha$ be the normal bundle of $\check{X}_\alpha$ as described in
(\ref{eq6.2}), a real line bundle, and consider the vector
bundle $\check{E}_\alpha^- \: = \: \check{E}_\alpha \otimes \nu_\alpha$
on $\check{X}_\alpha$. It is $G$-invariant and descends to a vector bundle
$E^-_\alpha$ on $X_\alpha$.

Although it is not strictly necessary for this paper, for explicitness
we wish to describe the elements of $\Omega(X; E)$ in more conventional terms 
on the stratified space 
$X$. If $\alpha < \alpha^\prime$ then there is a fiber bundle
$\pi_{\alpha \alpha^\prime} : F_{\alpha \alpha^\prime} 
\rightarrow X_\alpha$ such that for a small
tubular neighborhood $N_\alpha(\epsilon)$ of $X_\alpha$ in
$X$, $X_\alpha \cup (X_{\alpha^\prime} \cap N_\alpha(\epsilon))$ 
is homeomorphic to
the mapping cylinder of $\pi_{\alpha \alpha^\prime} : F_{\alpha \alpha^\prime} 
\rightarrow X_\alpha$. To see this, 
let $V(\epsilon)$ denote the $\epsilon$-ball in $V$. From (\ref{eq6.2}), a
tubular neighborhood $\check{N}_\alpha(\epsilon)$ of $\check{X}_\alpha$ is
$G$-diffeomorphic to
\begin{equation} \label{eq6.21}
\nu_\alpha(\epsilon) \: = \: P_\alpha \times_{S_\alpha} (V(\epsilon)
\times_H G).
\end{equation}
Let $(V \times_H G)_{\alpha^\prime}$ be the set of points in
$V \times_H G$ with normal orbit type $\alpha^\prime$ and put
$(V \times_H G)_{\alpha^\prime}^1 =
((V \times_H G)_{\alpha^\prime}) \cap (SV \times_H G)$. Then we can take
$F_{\alpha \alpha^\prime} = P_\alpha \times_{S_\alpha} 
(V \times_H G)_{\alpha^\prime}^1 /G$.

The vector bundle $E_{\alpha^\prime}$, when restricted to
$X_{\alpha^\prime} \cap N_\alpha(\epsilon)$,
extends to a vector bundle 
on $[0, \epsilon) \times F_{\alpha \alpha^\prime}$,
where we think of $\{0\} \times F_{\alpha \alpha^\prime}$ as the part which
is collapsed to $X_\alpha$ in the above mapping cylinder description. 
We write $E_{\alpha^\prime} \big|_{F_{\alpha \alpha^\prime}}$ for the
restriction of $E_{\alpha^\prime}$ to  $\{0\} \times 
F_{\alpha \alpha^\prime}$.

\begin{lemma} \label{lemmabefore}
There is an injection $I_{\alpha \alpha^\prime} :
\pi_{\alpha \alpha^\prime}^* E_\alpha \rightarrow 
E_{\alpha^\prime} \big|_{F_{\alpha \alpha^\prime}}$.
\end{lemma}
\begin{pf}
Using (\ref{eq6.2}),
it is enough to consider the case when $\check{X} = V \times_H G$
and $\check{E}$ is a $G$-equivariant vector bundle on $\check{X}$.
Then there is an $H$-module $W$ so that $\check{E} \: = \: 
(W \times V) \times_H G$. As $X \: = \:  \check{X}/G
\: = \: V/H$, for the purposes of the proof 
we may assume that $G \: = \: H$ and $\check{X} \: = \: V$. Taking
$\alpha \: = \: [H, V]$, we have that
$\check{X}_\alpha = 0 \in V$ and $X_\alpha$ is the point
$0 \cdot H$ in $X = V/H$. Suppose that 
$v \in SV \cap \check{X}_{\alpha^\prime}$ has isotropy group
$H^\prime \subset H$. Then the fiber of 
$E_{\alpha^\prime}$ over $vH \in X_{\alpha^\prime}$ is isomorphic to
$\check{K}_v^{H^\prime}$. The lemma follows from the 
injection $\check{K}_0^{H} \rightarrow \check{K}_0^{H^\prime}$,
along with our assumption that $\{\check{K}_{rv}\}_{r > 0}$ extends
continuously to $r \: = \: 0$.
\end{pf}

In the same way, if $\alpha \: = \: [H, \R]$ and
$\alpha^\prime \: = \: [H_\beta, 0]$ then there is an injection
$I_{\alpha \alpha^\prime} \: : \: \pi_{\alpha \alpha^\prime}^* 
E^-_\alpha \rightarrow
E_{\alpha^\prime} \big|_{F_{\alpha \alpha^\prime}}$.

Consider a collection of forms $\{\omega_\alpha \in
\Omega(X_\alpha; E_\alpha) \}_{\alpha \in {\cal N}}$.
If $\alpha < \alpha^\prime$ and 
$r \in [0, \epsilon)$ is the coordinate in the above mapping cylinder
description then we can write
$\omega_{\alpha^\prime}$ on $X_{\alpha^\prime} \cap N_\alpha(\epsilon)$ as
\begin{equation} \label{eq6.27}
\omega_{\alpha^\prime} \: = \: \omega_1(r) \: + \: dr \: \wedge \:
\omega_2(r),
\end{equation}
where for $r > 0$, we have 
$\omega_1(r), \omega_2(r) \in \Omega \left(
F_{\alpha \alpha^\prime}; E_{\alpha^\prime} 
\big|_{F_{\alpha \alpha^\prime}} \right)$.

We define a space of forms
$\Omega_{str}(X; E)$ which we call the stratified forms.
\begin{definition} \label{def8}
The forms $\{\omega_\alpha \in \Omega(X_\alpha; E_\alpha)
\}_{\alpha \in {\cal N}}$    
define an element of $\Omega_{str}(X; E)$ if
for $\alpha < \alpha^\prime$, the forms
$\omega_1(r)$ and $\omega_2(r)$ of (\ref{eq6.27})
are smooth up to $r = 0$, $\omega_1(0) = I_{\alpha \alpha^\prime}
\left( \pi_{\alpha \alpha^\prime}^* \omega_\alpha \right)$
and $\omega_2(0) \in
\begin{cases}
0 &\text{ if $\alpha \ne [H, \R]$,} \\
I_{\alpha \alpha^\prime}(\pi_{\alpha \alpha^\prime}^* 
E^-_\alpha) &\text{ if $\alpha = [H, \R]$.}
\end{cases}
$ 
\end{definition}
\noindent
{\bf Example 1 : } Let $X$ be a compact manifold-with-boundary. Let
$\check{X}$ be the double of $X$, with $G \: = \: \Z_2$ 
acting on $\check{X}$ so that the
generator $\gamma \in \Z_2$ acts by involution.  Put $\check{E} \: = \:
\check{X} \times \R$. If $\gamma$ acts on $\check{E}$ by sending
$(\check{x}, t)$ to $(\check{x} \gamma, t)$ then 
$E$ is the
trivial $\R$-bundle on $X$ and $\Omega_{str}(X; E)$
consists of the smooth forms on $X$ which lie in the quadratic
form domain of the
differential form Laplacian on $X$ with absolute boundary
conditions. If $\gamma$ acts on $\check{E}$ by sending
$(\check{x}, t)$ to $(\check{x} \gamma, - t)$ then the fiber of
$E$ is $\R$ over $X - \partial X$ and $0$ over $\partial X$, 
$E_\alpha^-$ is the trivial $\R$-bundle on $\partial X$ and 
$\Omega_{str}(X; E)$
consists of the smooth forms on $X$ which lie in the quadratic
form domain of the
differential form Laplacian on $X$ with relative boundary
conditions.\\ \\
{\bf Example 2 : } Let $Y$ be a compact manifold-with-boundary.
Take $G = S^1$ and $\check{X} = (Y \times S^1) \cup_{\partial Y \times S^1} 
(\partial Y \times D^2)$. Then the quotient space is the stratified space $X = 
\check{X}/S^1 = Y$. Let
$\check{E}$ be the trivial $\R$-bundle on $\check{X}$, with
${\frak I}_{\frak x} \equiv 0$. Then $E$ is the
trivial $\R$-bundle on $X$. One finds that $\Omega_{str}(X; E)$
consists of the smooth forms on $Y$ which lie in the quadratic
form domain of the
differential form Laplacian on $Y$ with absolute boundary
conditions.\\ \\
{\bf Example 3 : } Put $\check{X} = \C^2$ and $G = \Z_p \subset S^1$, with
the standard action on $\C^2$. (This is a noncompact example, but it 
illustrates the point.) Then $X = \cone(S^3/\Z_p)$. 
Let $\rho : \Z_p \rightarrow O(N)$ be a representation with
$(\R^N)^{\Z_p} = \{0\}$. Put $\check{E} = \C^2 \times \R^N$, with
the diagonal $\Z_p$-action and the
trivial product connection. Then
$E$ vanishes when restricted to the cone point of $X$, and has fiber
$\R^N$ on the rest of $X$. Putting ${\cal W} = S^3 \times_{\Z_p} \R^N$,
a flat vector bundle on $S^3/\Z_p$, one finds that $\Omega_{str}(X; E)$
consists of the elements of $\Omega^*([0, \infty)) \: \widehat{\otimes} \: 
\Omega^*(S^3/\Z_p; {\cal W})$ which vanish at $0$. That is, the entire
form vanishes at $0$ and not just its pullback.\\

\begin{lemma}
There is an inclusion $\Omega(X; E) \subset \Omega_{str}(X; E)$.
\end{lemma}
\begin{pf}
As in the proof of Lemma \ref{lemmabefore},
we can reduce to the case $\check{X} \: = \: V \: = \: V \times_H H$ and
$\check{E} \: = \: W \times V$, with $\check{X}_\alpha \: = \: 0 \in V$.
If $\check{\omega} \in
\Omega_{basic}(\check{X}; \check{E})$ then write
$\check{\omega} \: = \: \check{\omega}_1(r) \: + \: dr \: \wedge \:
\check{\omega}_2(r)$, where $r$ is the radial coordinate on $V$.

Suppose that $\alpha \ne [H, \R]$.
As $\check{\omega}$ is smooth and $\dim(V) \: \ge \: 2$, 
we must have that $\check{\omega}_2(0) \: = \: 0$
and $\check{\omega}_1(0) \in \Lambda^0(T^*V) \otimes W$.
Then as $\check{\omega}$ is basic,
we have $\check{\omega}_1(0) \in \Lambda^0(T^*V) \otimes \check{K}_0^H$.
Suppose that $v \in SV \cap \check{X}_{\alpha^\prime}$ 
has isotropy group $H^\prime \subset H$.
By the smoothness of
$\check{\omega}$, we have
\begin{equation} \label{trans}
\lim_{r \rightarrow 0} \check{\omega}_1(rv) \: = \: \check{\omega}_1(0).
\end{equation}
In particular, $\lim_{r \rightarrow 0} \check{\omega}_1(rv)$ is a $0$-form
with value in $\check{K}_0^H \subset \check{K}_0^{H^\prime}$.
Translating (\ref{trans}) to the quotient space gives
$\omega_1(0) = I_{\alpha \alpha^\prime}
\left( \pi_{\alpha \alpha^\prime}^* \omega_\alpha \right)$.

If $\alpha = [H, \R]$ then the only difference is that
$\check{\omega}_2(0)$ can be nonzero.  From $H$-invariance, it
must lie in $\check{E}^-_\alpha$.
\end{pf}

Hence we have a space of forms $\Omega_{str}(X; E)$ which is defined
intrinsically on the stratified space $X$ and which in some sense is the
minimal such space that contains the basic forms $\Omega(X; E)$.

Finally, we describe a spectral sequence to compute the cohomology of
$\check{A}^\prime$ acting on $\Omega_{basic}(\check{X}; \check{E})$,
which we denote by $\HH^*(A^\prime)$.
In the case when $G \: = \: \{e\}$, $\check{X} \: = \: B$ and 
$\check{E} \: = \: E$, 
such a spectral sequence was described in
\cite[Section 7]{Lott (2001)}. It arises from the filtration of
$\Omega \left( B; E \right)$ 
by $F^p \: = \: \bigoplus_{q \ge p} \Omega^q(B; E^*)$.
In our case, we
filter $\Omega_{basic}(\check{X}; \check{E})$ by 
the form degree on the quotient space $X$. More precisely,
let $F^p$ be the forms $\check{\omega} \in
\Omega_{basic}(\check{X}; \check{E})$ such that for all
$\alpha \in {\cal N}$, $\omega_\alpha \in 
\bigoplus_{q \ge p} \Omega^q(X_\alpha; E_\alpha)$. As usual, the
$E_1$-term $E_1^{p,q}$ of the spectral sequence is the $(p \: + \: q)$-degree
cohomology of the complex $F^p/F^{p+1}$. Given an open set $U \subset X$,
let $\check{U}$ be its preimage in $\check{X}$. Let $E_1^{p,q}$ also denote
the sheaf which assigns to $U$ the 
space $E_1^{p,q}(U)$ as computed using the basic forms on $\check{U}$.
If $\check{x} \in \check{X}$ covers $x \in X$ and has isotropy group $H$ then
the stalk of $E_1^{0,*}$ at $x$ is isomorphic to the cohomology of
$\check{A}^\prime_{[0]}$ on $\check{K}_{\check{x}}^H$. 
(The degree-$1$ component of the equation
$\check{A}^\prime \: (i_{\frak X} \: + \: {\frak I}_{\frak x}) \: + \:
(i_{\frak X} \: + \: {\frak I}_{\frak x}) \: \check{A}^\prime \: = \:
{\cal L}_{\frak X}$ is
$\check{A}^\prime_{[1]} \: i_{\frak X} \: + \:
i_{\frak X} \: \check{A}^\prime_{[1]} \: + \:
{\frak I}_{\frak x} \: \check{A}^\prime_{[0]} \: + 
\: \check{A}^\prime_{[0]} \: {\frak I}_{\frak x}  \: = \:
{\cal L}_{\frak X}$. 
Taking ${\frak x} \in {\frak h}$, one sees that
$\check{A}^\prime_{[0]}$ does act on $\check{K}_{\check{x}}^H$.)

Define a $\Z$-graded sheaf $\HH^*(A^\prime_{[0]})$ on $X$ which assigns
to $U \subset X$ 
the vector space $\Ker (d_1^{0,*} \: : \: E_1^{0,*}(U) \rightarrow 
E_1^{1,*}(U))$. There is a complex of sheaves
\begin{equation} \label{complex}
\HH^*(A^\prime_{[0]}) \longrightarrow E_1^{0,*}
\stackrel{d_1^{0,*}}{\longrightarrow} 
E_1^{1,*} \stackrel{d_1^{1,*}}{\longrightarrow} E_1^{2,*}
\stackrel{d_1^{2,*}}{\longrightarrow} \ldots
\end{equation}
\begin{lemma} \label{Verona}
The complex (\ref{complex}) is a resolution of
$\HH^*(A^\prime_{[0]})$ by fine sheaves.
\end{lemma}
\begin{pf}
We follow the method of proof of
\cite{Verona (1988)}, which effectively proves the lemma in the
case when $\check{E}$ is the trivial $\R$-bundle on $\check{X}$. 
Using the ordinary slice theorem, we can reduce to the case
 that $\check{U} \: = \:
N \times_H G$ for some representation space $N$ of an isotropy group $H$.
Then we can reduce to the case when $\check{U} \: = \: N$ and
$G \: = \: H$. Now
$E_1^{p,*}(U)$ is the cohomology of $\check{A}^\prime_{[0]}$ on the elements
of $\Omega_{basic} \left( N; \check{E} \right)$ with $p$ horizontal
differentials, with respect to the $H$-action on $N$. The degree-$1$
component of the equation
$ \left( \check{A}^\prime \right)^2
\: = \: 0$ is $\check{A}^\prime_{[1]} \: \check{A}^\prime_{[0]} \: + \: 
\check{A}^\prime_{[0]} \: \check{A}^\prime_{[1]} \: = \: 0$, which implies
that $\check{A}^\prime_{[1]}$ has an induced action on $E_1^{p,*}(U)$. 
The degree-$2$ component of the equation $ \left( \check{A}^\prime \right)^2
\: = \: 0$ is $\left( \check{A}^\prime_{[1]} \right)^2 \: + \: 
\check{A}^\prime_{[0]} \: \check{A}^\prime_{[2]} \: + \: 
\check{A}^\prime_{[2]} \: \check{A}^\prime_{[0]} \: = \: 0$, which implies
that $\left( \check{A}^\prime_{[1]} \right)^2$ vanishes on $E_1^{p,*}(U)$.
In fact, the action of $\check{A}^\prime_{[1]}$ on 
$E_1^{p,*}(U)$ is the same as $d_1^{p,*}$. We now use the
Poincar\'e lemma on $N$, as in \cite{Verona (1988)}, to prove the
claim.  To apply the Poincar\'e lemma we use a radial trivialization of
$\check{A}^\prime_{[1]}$ on $N$. Thinking of $N$ as the cone over its
sphere $SN$, the $H$-action on $N$ comes from the $H$-action on
$SN$. Because of this, it follows that the homotopy operator in the
Poincar\'e lemma does send basic forms to basic forms. 
The rest of the proof is as in \cite{Verona (1988)}.
\end{pf}

It follows that
the $E_2$-term of the spectral sequence is given by
$E_2^{p,q} \: = \: \HH^p \left( X; \HH^q(A^\prime_{[0]}) \right)$, where
the right-hand-side is the cohomology of the sheaf
$\HH^q(A^\prime_{[0]})$ on $X$.

\section{Eigenvalue Bounds} \label{sect7}

In this section we use the results of Section \ref{sect6} to 
prove the analogs of \cite[Theorems 2 and 3]{Lott (2001)}
in the case of a general limit space $X$.
We then prove Theorem \ref{th5}, giving eigenvalue
bounds for the $p$-form Laplacian in terms of sectional curvature
and diameter. The method of proof is to assume that there are no such bounds
and use Gromov-Hausdorff convergence in the $O(n)$-equivariant setting,
along with our eigenvalue estimates, to get a contradiction.
In Theorem \ref{tcor1} we look at the special case of Theorem
\ref{th5}.2 when $p = 2$.

Let $M$ be a smooth connected closed $n$-dimensional
Riemannian manifold and let $FM$ denote
its orthonormal
frame bundle, the total space of a principal $O(n)$-bundle
${\frak p} \: : \: FM \rightarrow M$.
There is a canonical Riemannian metric on $FM$, but for the moment
we will consider $FM$ with any $O(n)$-invariant Riemannian metric, and
we give $M$ the corresponding quotient metric.
Let $\mu$ be the (smooth) measure on $M$ given by
$\mu(m) \: = \: \vol \left({\frak p}^{-1}(m) \right) \: d\vol(m)$.
Let $\Omega^*_{L^2}(M, \mu)$ denote the completion of $\Omega^*(M)$ with
respect to the inner product
$\langle \omega, \omega \rangle \: = \: 
\int_M |\omega(m)|^2 \: d\mu(m)$.
Then there is an isometric isomorphism
$\Omega^*_{L^2}(M, \mu) \: \cong \:  \Omega^*_{basic,L^2}(FM)$ coming from
pullback.   
The basic Laplacian $\triangle^{FM}_{basic}$ 
on $\Omega^*_{basic, L^2}(FM)$ is the Laplacian associated to the
complex
\begin{equation} \label{complexx}
\Omega^0_{basic}(FM) \longrightarrow \Omega^1_{basic}(FM) \longrightarrow
\ldots
\end{equation}
It is isomorphic to the 
weighted Laplacian $\triangle^M_{\mu} \: = \: d d^* \: + \: d^* d$ on 
$\Omega^*_{L^2}(M, \mu)$.

Let $\check{X}$ be a fixed
smooth connected closed Riemannian manifold on which $O(n)$ acts 
isometrically on the right, with quotient $X = \check{X}/O(n)$. 
We say that a fiber bundle $FM \rightarrow
\check{X}$ is an equivariant Riemannian affine fiber bundle if
it is a Riemannian affine fiber bundle
\cite[Definition 1]{Lott (2001)} which is $O(n)$-equivariant in
the obvious sense.
Given $\check{x} \in \check{X}$, let
$\check{Z}_{\check{x}}$ be the fiber over $\check{x}$ of the affine
fiber bundle. For the applications, it will be sufficient to consider
the case when $\check{Z}_{\check{x}}$ is a nilmanifold
$\Gamma \backslash N$ \cite[(7.2)]{Cheeger-Fukaya-Gromov (1992)}.
Let $\check{E}$ be the real
$\Z$-graded vector bundle on $\check{X}$ whose fibers are isomorphic
to the affine-parallel differential
forms on $\{ \check{Z}_{\check{x}} \}_{\check{x} \in
\check{X}}$. It inherits a flat degree-$1$ superconnection
$\check{A}^\prime$ from $d^{FM}$ \cite[Section 5]{Lott (2001)}.
If ${\frak x}$ is in the Lie algebra $o(n)$ of $O(n)$, let ${\frak X}^{FM}$ be
the corresponding vector field on $FM$ and let ${\frak X}_V$ be its 
vertical component with respect to $FM \rightarrow \check{X}$. Then
interior multiplication by ${\frak X}_V$ on the fibers $\check{Z}_{\check{x}}$
induces a linear map ${\frak I}_{\frak x} \in 
C^\infty \left( \check{X}; \Hom ( \check{E}^*, 
\check{E}^{*-1}) \right)$, which gives
$\check{E}$ the structure of an $O(n)$-basic vector bundle.
Furthermore, $\check{A}^\prime$ becomes an $O(n)$-basic superconnection.
Define $\triangle^{\check{E}}_{basic} \cong \triangle^E$ as in Definition 
\ref{def7}. It is the Laplacian associated to the subcomplex
$\Omega_{basic} \left( \check{X}; \check{E} \right)$ of (\ref{complexx}).

The Riemannian metric on $FM$ defines an $O(n)$-invariant
family of horizontal planes
that are perpendicular to the fibers of $FM \rightarrow M$.
Let $T^{FM}$ be the corresponding fiber bundle curvature.
Let $\check{T}$ be the curvature of the affine fiber bundle
$FM \rightarrow \check{X}$, let $\check{\Pi}$ be
the second fundamental forms of the fibers 
$\{ \check{Z}_{\check{x}} \}_{\check{x} \in
\check{X}}$ and let $\diam(\check{Z})$ be the maximum
diameter of the fibers.

The next proposition can be considered to be an analog of
\cite[Theorem 1]{Lott (2001)}, in which the nilpotent fiber bundle structure
on the total space $M$ of an affine fiber bundle $M \rightarrow B$ is replaced
by the nilpotent Killing structure on $M$ coming from an $O(n)$-equivariant
affine fiber bundle $FM \rightarrow \check{X}$
\cite{Cheeger-Fukaya-Gromov (1992)}.

\begin{proposition} \label{prop10}
There are positive constants $A$, $A^\prime$ and 
$C$ which only depend on $\dim(M)$ such that
if $\parallel R^{\check{Z}} \parallel_\infty 
\diam( \check{Z})^2 \le A^\prime$ then
\begin{align} \label{eq7.1}
& \sigma(\triangle^{FM}_{basic}) \cap \left[ 0, 
\left( 
\sqrt{ A \: 
{\diam(\check{Z})^{-2}} \: - \: C \: \left(
\parallel R^{FM} \parallel_\infty + 
\parallel \check{\Pi} \parallel_\infty^2 +
\parallel \check{T} \parallel_\infty^2 
 \right) } \: - \: C \: \parallel T^{FM} \parallel_\infty \right)^2 \right)
= \\
& \sigma(\triangle^E) \cap \left[ 0, 
\left( 
\sqrt{ A \: 
{\diam(\check{Z})^{-2}} \: - \: C \: \left(
\parallel R^{FM} \parallel_\infty + 
\parallel \check{\Pi} \parallel_\infty^2 +
\parallel \check{T} \parallel_\infty^2 
 \right) } \: - \: C \: \parallel T^{FM} \parallel_\infty \right)^2 \right).
\notag
\end{align}
\end{proposition}
\begin{pf}
Let ${\cal P}^{fib}$ be fiberwise orthogonal projection from
$\Omega^*(FM)$ to $\Omega(\check{X}; \check{E})$.
From
the proof of \cite[Proposition 1]{Lott (2001)},
${\cal P}^{fib}$ amounts to averaging
over the nilmanifold fibers $\Gamma \backslash N$
of $FM \rightarrow \check{X}$. Hence it 
preserves $\Omega^*_{basic}(FM)$ and descends
to an operator from $\Omega^*_{basic}(FM)$ to 
$\Omega_{basic}(\check{X}; \check{E}) \cong 
\Omega(X; E)$, which we also denote by ${\cal P}^{fib}$. As in the proof
of \cite[Theorem 1]{Lott (2001)}, it suffices to show that 
there are constants $A$, $A^\prime$ and $C$ as above such that 
$\sigma \left( \triangle^{FM}_{basic} 
\big|_{Ker({\cal P}^{fib})} \right)$ is bounded below by
\begin{equation} \label{eq7.2}
\left( 
\sqrt{ A \: 
{\diam(\check{Z})^{-2}} \: - \: C \: \left(
\parallel R^{FM} \parallel_\infty + 
\parallel \check{\Pi} \parallel_\infty^2 +
\parallel \check{T} \parallel_\infty^2 
 \right) } \: - \: C \: \parallel T^{FM} \parallel_\infty \right)^2.
\end{equation}

Let $\triangle^{FM}$ be the Laplacian on differential forms on
$FM$ and let $\triangle^{FM}_{O(n)}$ be the 
Laplacian on $O(n)$-invariant not-necessarily-basic differential forms on
$FM$. Applying \cite[Theorem 1]{Lott (2001)} to the Riemannian
affine fiber bundle $FM \rightarrow
\check{X}$, we know that
\begin{equation} \label{eq7.3}
\sigma \left( \triangle^{FM}
\big|_{Ker({\cal P}^{fib})} \right) \: \subset \: 
\left[
A \: 
{\diam(\check{Z})^{-2}} \: - \: C \: \left(
\parallel R^{FM} \parallel_\infty + 
\parallel \check{\Pi} \parallel_\infty^2 +
\parallel \check{T} \parallel_\infty^2 
 \right), \infty \right) 
\end{equation}
and so
\begin{equation} \label{eq7.4}
\sigma \left( \triangle^{FM}_{O(n)}
\big|_{Ker({\cal P}^{fib})} \right) \: \subset \: 
\left[
A \: 
{\diam(\check{Z})^{-2}} \: - \: C \: \left(
\parallel R^{FM} \parallel_\infty + 
\parallel \check{\Pi} \parallel_\infty^2 +
\parallel \check{T} \parallel_\infty^2 
 \right), \infty \right) 
\end{equation}

Let $d^*_{O(n)}$ denote the adjoint of $d$ on $\Omega^*_{O(n)}(FM)$. 
As $\triangle^{FM}_{O(n)} \: = \: \left( d + d^*_{O(n)} \right)^2$,
it follows from (\ref{eq7.4}) that if $\lambda$ is an eigenvalue of 
$\left( d + d^*_{O(n)} \right)\big|_{Ker({\cal P}^{fib})}$ then 
\begin{equation} \label{eq7.5}
|\lambda| \: \ge \: \sqrt{
A \: 
{\diam(\check{Z})^{-2}} \: - \: C \: \left(
\parallel R^{FM} \parallel_\infty + 
\parallel \check{\Pi} \parallel_\infty^2 +
\parallel \check{T} \parallel_\infty^2 
 \right)
}.
\end{equation}

From Lemma \ref{lemma5}, the adjoint of $d$ on 
$\Omega^*_{basic}(FM) \cong \Omega^*(M)$ is
$d^*_{basic} \: = \: P^{hor} \: d^*_{O(n)}$. 
With respect to the decomposition
$\Omega^*_{O(n)}(FM) \: = \: \Omega^*_{basic}(FM) \: \oplus \:
\left( \Omega^*_{basic}(FM) \right)^\perp$, we have
\begin{equation} \label{eq7.6}
d \: + \: d^*_{O(n)} \: = \:
\begin{pmatrix}
d \: + \: d^*_{basic} & P^{hor} \: d \: ( I \: - \: P^{hor}) \\
( I \: - \: P^{hor}) \:  d^*_{O(n)} P^{hor} & 
( I \: - \: P^{hor}) \: \left( d \: + \: d^*_{O(n)}
\right) \: ( I \: - \: P^{hor})
\end{pmatrix}.
\end{equation}
Using the notation of \cite[(5.26)]{Lott (2001)}, let 
$\omega_{i \alpha \beta}$ denote
the curvature of the fiber bundle $FM \rightarrow M$. A calculation gives
that 
\begin{equation} \label{eq7.7}
( I \: - \: P^{hor}) 
\: d^*_{O(n)} \: P^{hor} \: = \: - \: \sum_{i,\alpha,\beta}
\omega_{i \alpha \beta} \: E^i \: I^\alpha \: I^\beta,
\end{equation}
with $P^{hor} \: d \: ( I \: - \: P^{hor})$ being the adjoint of the
right-hand-side of (\ref{eq7.7}).
Then upon restriction to
$\Ker({\cal P}^{fib})$, it follows that there for an appropriate 
constant $C = C(\dim(M)) > 0$,
the diagonal part of the operator in
(\ref{eq7.6}) has a spectrum which differs from that of
$\left( d + d^*_{O(n)} \right)\big|_{Ker({\cal P}^{fib})}$ by at most
$C \: \parallel T^{FM} \parallel_\infty$. As the spectrum of
$\left( d + d^*_{basic} \right)\big|_{Ker({\cal P}^{fib})}$ is contained in the
spectrum of the diagonal part of (\ref{eq7.6}), when restricted to
$\Ker({\cal P}^{fib})$, the proposition follows.
\end{pf}

Let $d^{O(n)}_{GH}$ denote the $O(n)$-equivariant
Gromov-Hausdorff metric on the space of $O(n)$-equivariant
compact metric spaces 
\cite[(2.1.3)]{Cheeger-Fukaya-Gromov (1992)}. We say that two nonnegative
numbers $\lambda_1$, $\lambda_2$ are $\epsilon$-close if
$e^{- \: \epsilon} \: \lambda_2 \: \le \: \lambda_1 \: \le \:
e^{\epsilon} \: \lambda_2$.

\begin{proposition} \label{prop11} 
For $n \in \Z^+$, let $\check{X}$ be a fixed
connected closed
Riemannian manifold with an isometric $O(n)$-action and
quotient space $X = \check{X}/O(n)$.  
Given $\epsilon > 0$ and
$K \ge 0$, there are positive constants $A(n, \epsilon, K)$,
$A^\prime(n, \epsilon, K)$ and 
$C(n, \epsilon, K)$ with the following property :
If $M^n$ is an $n$-dimensional connected closed
Riemannian manifold with $\parallel R^M \parallel_\infty
\: \le \: K$ and $d^{O(n)}_{GH}(FM, \check{X}) \: \le 
\: A^\prime(n,\epsilon, K)$
then there are\\
1. An $O(n)$-basic $\Z$-graded real vector bundle 
$\check{E}$ on $\check{X}$,\\
2. A basic flat degree-$1$ superconnection 
$\check{A}^\prime$ on $\check{E}$ and\\
3. An $O(n)$-invariant Euclidean inner product $h^{\check{E}}$ on 
$\check{E}$\\
such that    
if $\lambda_{p,j}(M)$ is the $j$-th eigenvalue of the $p$-form Laplacian
on $M$, $\lambda_{p,j}(X; E)$ is the $j$-th eigenvalue of
$\triangle^E_p$ and
\begin{equation} \label{eq7.8}
\min (\lambda_{p,j}(M), \lambda_{p,j}(X; E)) \: \le 
\: \left( 
\sqrt{ A(n, \epsilon, K) \: 
{d_{GH}^{O(n)}(FM, \check{X})^{-2}} \: - \: C(n, \epsilon, K) } 
\: - \: C(n, \epsilon, K) \right)^2 
\end{equation}
then $\lambda_{p,j}(M)$ 
is $\epsilon$-close to $\lambda_{p,j}(X; E)$.
\end{proposition}
\begin{pf}
For simplicity, we will omit reference to $p$. Let $g_0^{M}$ denote
the Riemannian metric on $M$. Let $g_0^{FM}$ denote
the canonical induced Riemannian metric on $FM$.
Then $\triangle^M \: \cong \: \triangle^{FM}_{basic}$. 

From the obvious basic extension of \cite[Lemma 3]{Lott (2001)},
if an $O(n)$-invariant Riemannian metric $g_1^{FM}$ on $FM$ is $\epsilon$-close
to $g_0^{FM}$ in the sense of \cite[(5.4)]{Lott (2001)} 
then the spectrum of $\triangle^{FM}_{basic}$ 
computed with $g_1^{FM}$
is $J\epsilon$-close to the spectrum computed with $g_0^{FM}$. 
We can use the geometric results of \cite{Cheeger-Fukaya-Gromov (1992)} to
find a $O(n)$-invariant Riemannian metric $g_2^{FM}$ on $M$ which is close to 
$g_0^{FM}$ and to which we can apply Proposition \ref{prop10}. Note that
$g_2^{FM}$ may not be the canonical metric coming from some Riemannian
metric on $M$.

There is an explicit bound for the sectional curvatures of the canonical
metric $g_0^{FM}$ in terms of $K$. Then in the construction of the 
$O(n)$-invariant metric $g_2^{FM}$, we may assume that we have bounds
on the sectional curvatures of $g_2^{FM}$ and on the sectional
curvatures of the quotient metric on $M$, in terms of $\epsilon$ and $K$
(see \cite[Proof of Theorem 1.3]{Cheeger-Fukaya-Gromov (1992)} and
\cite[Theorem 2.1]{Rong (1996)}).

We can go through the proof of \cite[Theorem 2]{Lott (2001)} 
working $O(n)$-equivariantly
on $FM$ and using Proposition \ref{prop10}. The proofs of 
\cite[Propositions 3.6 and 4.9]{Cheeger-Fukaya-Gromov (1992)}, giving
the equivariant Riemannian affine fiber bundle structure $FM \rightarrow
\check{X}$, are explicitly phrased in the $G$-equivariant setting,
where $G$ is a compact Lie group. In particular, we have bounds on
$\parallel \check{\Pi} \parallel_\infty$ and 
$\parallel \check{T} \parallel_\infty$. 
Using O'Neill's formula \cite[(9.29c)]{Besse (1987)}, we have a bound
on $\parallel {T}^{FM} \parallel_\infty$ in terms of
$\epsilon$ and $K$.
Then the proof of 
\cite[Theorem 2]{Lott (2001)} goes through to the $O(n)$-equivariant
setting, changing the
fibration $M \rightarrow B$ of \cite[Theorem 2]{Lott (2001)}
to the $O(n)$-equivariant fibration
$FM \rightarrow \check{X}$.
\end{pf}

If $\check{E}$ is a 
real $\Z$-graded $O(n)$-equivariant
vector bundle on $\check{X}$, let
${\cal I}_{\check{E},basic}$ be the set of $O(n)$-equivariant linear maps
${\frak I} \: : \: o(n) \rightarrow 
C^\infty \left( \check{X}, \Hom(\check{E}^*, \check{E}^{*-1}) \right)$
which satisfy ${\frak I}({\frak x})^2 \: = \: 0$ for all
${\frak x} \in o(n)$.
Let ${\cal S}_{\check{E}, O(n)}$ be the set of smooth degree-$1$
$O(n)$-equivariant superconnections on $\check{E}$ and
let ${\cal S}_{\check{E},basic}$ be the set of pairs
$({\frak I}, \check{A}^\prime) \in {\cal I}_{\check{E},basic} \times
{\cal S}_{\check{E}, O(n)}$ so that $\check{A}^\prime$ is $O(n)$-basic
with respect to ${\frak I}$.
Let ${\cal H}_{\check{E},basic}$ be the set of ${O(n)}$-invariant
Euclidean inner products on 
$\check{E}$ and
let ${\cal G}_{\check{E},basic}$ be the group of smooth 
grading-preserving ${O(n)}$-equivariant
$\GL(\check{E})$-gauge transformations on $\check{E}$ and
We equip ${\cal I}_{\check{E},basic}$, 
${\cal S}_{\check{E},O(n)}$ and
${\cal H}_{\check{E},basic}$
with the $C^\infty$-topology, ${\cal S}_{\check{E},basic}$ with the
subspace topology and
$( {\cal S}_{\check{E},basic} \times
{\cal H}_{\check{E},basic} )/ {\cal G}_{\check{E},basic}$ with the
quotient topology.

\begin{proposition} \label{prop12}
In Proposition \ref{prop11}, we may assume that $\check{E}$ 
is in one of a finite number of isomorphism classes of
real $\Z$-graded $O(n)$-equivariant topological vector bundles 
$\{\check{E}_i\}$ on $\check{X}$. In addition, we may assume that for each
$\alpha \in {\cal N}$,
the vector bundle $E_\alpha$ is in one of a finite number of isomorphism
classes of real $\Z$-graded topological vector bundles on $X_\alpha$.
Furthermore, there are compact subsets $C_{\check{E}_i} \subset
({\cal S}_{\check{E}_i,basic} \times 
{\cal H}_{\check{E}_i,basic})/{\cal G}_{\check{E}_i,basic}$ 
depending on $n$, $\epsilon$ and $K$ such that we 
may assume that the gauge-equivalence class of 
$\left({\frak I}, \check{A}^\prime, h^{\check{E}} \right)$
lies in $C_{\check{E}}$.
\end{proposition}
\begin{pf}
We go through the proof of 
\cite[Theorem 3]{Lott (2001)} equivariantly.
In \cite[Theorem 3]{Lott (2001)}, one obtains the finiteness statement
from the fact that there is a finite number of topological types of
real vector bundles of a given rank on $B$ which admit a flat connection.
This fact follows from the finiteness of the number of connected 
components of the representation variety of $\pi_1(B, b_0)$. In
our case the representation variety of $\pi_0(\widehat{G})$ will have
a finite number of connected components. With Proposition
\ref{prop6}, this implies that there is a finite number of topological
types of $O(n)$-equivariant vector bundles on $\check{X}$ which admit
a basic flat connection. The method of proof of \cite[Theorem 3]{Lott (2001)}
shows that $\check{E}$ admits an $O(n)$-basic flat connection.

The vector bundle $E_\alpha$ on $X_\alpha$ has a flat degree-$1$
superconnection $A^\prime_\alpha$ induced from $\check{A}^\prime$.  From
the argument of \cite[Theorem 3]{Lott (2001)}, there is a finite number
of possibilities for the topological type of $E_\alpha$.

Next, \cite[Theorem 3]{Lott (2001)} gives a compactness
result for $\left( \check{A}^\prime, h^{\check{E}} \right)$ in
$({\cal S}_{\check{E}_i,O(n)} \times 
{\cal H}_{\check{E}_i,basic})/{\cal G}_{\check{E}_i,basic}$.
We recall that the inner product $h^{\check{E}}$ of Proposition \ref{prop12},
when restricted to the fiber $\check{E}_{\check{x}}$ over a point
$\check{x} \in \check{X}$, is the $L^2$-inner product
on the affine-parallel forms of the geometric fiber $\check{Z}_{\check{x}}$.  
As the $O(n)$-orbits on $FM$, with the canonical metric,
are isometric to the standard $O(n)$,
from the construction of ${\frak I}_{\frak x}$ 
we have an upper bound on 
\begin{equation} \label{upper}
\parallel {\frak I} \parallel^2 \: = \:
\sup_{{\frak x} \in o(n) 
\: : \: |{\frak x}| = 1} \:
\sup_{\check{x} \in \check{X}} \:
\sup_{\check{e} \in \check{E}_{\check{x}} - 0} 
\frac{h^{\check{E}}({\frak I}_{\frak x} \check{e},{\frak I}_{\frak x} 
\check{e})}{h^{\check{E}}(\check{e}, \check{e})}
\end{equation}
that only depends on $n$, $\epsilon$ and $K$.
The component of the equation
$\check{A}^\prime \: ({\frak I}_{\frak x} \: + \: {\frak i}_{\frak X}) \: + \:
({\frak I}_{\frak x} \: + \: {\frak i}_{\frak X}) \: \check{A}^\prime \: = \:
{\cal L}_{\frak X}$
of degree $1$ with respect to $\check{X}$ is
\begin{equation}
\check{A}^\prime_{[1]} \: {\frak I}_{\frak x} \: + \:
{\frak I}_{\frak x} \: \check{A}^\prime_{[1]} \: +
\check{A}^\prime_{[2]} \: {\frak i}_{\frak X} \: + \: {\frak i}_{\frak X}
\: \check{A}^\prime_{[2]}
\: = \:  0.
\end{equation}
As we have $C^\infty$-bounds on $\check{A}^\prime_{[2]}$, we obtain
$C^\infty$-bounds on the covariant derivative of ${\frak I}$, and
hence on its higher covariant derivatives. Thus we also have precompactness for
${\frak I}$, from which the proposition follows.
\end{pf}

We will need an eigenvalue estimate. Let $\check{E}$ be a real
$\Z$-graded $O(n)$-equivariant vector bundle on $\check{X}$ and define
${\cal S}_{\check{E},basic}$ and ${\cal H}_{\check{E},basic}$ as before.
Suppose that we have two triples
$\left( {\frak I}_1, \check{A}^\prime_1, h^{\check{E}}_1 \right)$ and
$\left( {\frak I}_2, \check{A}^\prime_2, h^{\check{E}}_2 \right)$ 
in ${\cal S}_{\check{E},basic} \times {\cal H}_{\check{E},basic}$.
For $i \in \{1,2\}$,
let $\Omega_{basic,i} \left( \check{X}; \check{E} \right)$ denote the
basic forms as defined using ${\frak I}_i$. 
Suppose that $T \: : \: \Omega_{basic,1} \left( \check{X}; \check{E} \right)
\rightarrow
\Omega_{basic,2} \left( \check{X}; \check{E} \right)$ is an isomorphism.
Let $\triangle^E_i$ denote the basic Laplacian constructed using
$\left( {\frak I}_i, \check{A}^\prime_i, h^{\check{E}}_i \right)$.

Given
$y \in \Omega_{basic,2}( \check{X}; \End({\check{E}}))$, 
let $\parallel y \parallel$ be the 
operator norm for the action of $y$ on $\Omega_{basic,2,L^2}
(\check{X}; {\check{E}})$.

\begin{lemma} \label{lemma??}
Suppose that $\triangle^E_1$ has a discrete spectrum. 
If for some $\epsilon \: \ge \: 0$ we have
$e^{- \: \epsilon} \: \Id. \: \le \: T^* T \: \le \: 
e^{\epsilon} \: \Id.$ then 
for all $j \in \Z^+$,
\begin{equation} \label{eq7.9}
|\lambda_j(\triangle^E_1)^{1/2} - 
\lambda_j(\triangle^E_2)^{1/2}| \: \le \:
(2 + \sqrt{2}) \: \parallel T \: \check{A}_1^\prime \: T^{-1} \:  - \: 
\check{A}_2^\prime \parallel \: + \: (e^{\epsilon} \: - \: 1) \:
\lambda_j(\triangle^E_1)^{1/2}.
\end{equation}
\end{lemma}
\begin{pf}
We first examine the effect on the eigenvalues
of changing from $(\check{A}_2^\prime,
h^{\check{E}}_2)$ to $(T \check{A}_1^\prime T^{-1},
h^{\check{E}}_2)$, where in both cases the corresponding
superconnection Laplacian
acts on $\Omega_{basic,2} \left( \check{X}; \check{E} \right)$. 
From the method of proof of \cite[Lemma 4]{Lott (2001)}
the change in $\lambda_j^{1/2}$ is bounded by
$(2 + \sqrt{2}) \: \parallel T \: \check{A}_1^\prime \: T^{-1} \:  - \: 
\check{A}_2^\prime \parallel$. Next, we examine the effect of
changing from $(T \check{A}_1^\prime T^{-1},
h^{\check{E}}_2)$ to $(\check{A}_1^\prime,
h^{\check{E}}_1)$. By naturality,
the eigenvalues of the superconnection Laplacian
constructed from
$(T \check{A}_1^\prime T^{-1},
h^{\check{E}}_2)$ can be computed using instead
the superconnection $\check{A}_1$
and the inner product on
$\Omega_{basic,1} \left( \check{X}; \check{E} \right)$ which is pullbacked
from $\Omega_{basic,2} \left( \check{X}; \check{E} \right)$ via $T$.
The method of proof of \cite[Lemma 3]{Lott (2001)} shows that if one
compares this with the original inner product on
$\Omega_{basic,1} \left( \check{X}; \check{E} \right)$ then the eigenvalues
differ at most by a multiplicative factor of $e^{2 \epsilon}$.
The lemma follows. We note that we have implicitly shown that
$\triangle^E_2$ also has a discrete spectrum.
\end{pf}
\noindent
{\bf Proof of Theorem \ref{th5} :} \\
{\bf 1.} 
If it is not true that $\lim_{j \rightarrow \infty} a^1_{n,p,j,K}$ is always
infinite then there are numbers $n \in \Z^+$, $0 \le p \le n$, $K \ge 0$, 
$\Lambda > 0$ and a sequence
$\{M_i\}_{i = 1}^\infty$ of connected closed $n$-dimensional Riemannian
manifolds with $\diam(M_i) = 1$ and 
$\parallel R^{M_i} \parallel_\infty \: \le \: K$ such that
$\lambda_{p,i}(M_i) < \Lambda$.
By \cite[Theorem 1.12]{Cheeger-Fukaya-Gromov (1992)}, for any $\epsilon > 0$
there is a sequence $\{A_k(n, \epsilon)\}_{k=0}^\infty$ so that we can
find a new metric on $M_i$ which is $\epsilon$-close to the old one
in the sense of \cite[(5.4)]{Lott (2001)},
with the new metric satisfying 
$\parallel \nabla^k R^{M_i} \parallel_\infty \: \le \:
A_k(n, \epsilon)$. 
Fix $\epsilon$ to be, say, $\frac{1}{2}$ and consider $\{M_i\}_{i = 1}^\infty$
with the new metrics. From \cite{Dodziuk (1982)} or 
\cite[Lemma 3]{Lott (2001)}, 
we now have
$\lambda_{p,i}(M_i) \: < \: e^{J\epsilon} \: \Lambda$ for a fixed integer $J$.
As in 
\cite[III.5]{Cheeger-Fukaya-Gromov (1992)}, we can apply
Gromov's convergence theorem
in the equivariant setting to conclude that there are a smooth Riemannian
$O(n)$-manifold $\check{X}$ and a subsequence of
$\{M_i\}_{i = 1}^\infty$, which we relabel as $\{M_i\}_{i = 1}^\infty$, so
that $d_{GH}^{O(n)} (FM_i, \check{X}) \: \le \: \frac{1}{i}$. In particular,
$X = \check{X}/O(n)$ 
is not a point. As in the proof of Proposition \ref{prop11}, we slightly
perturb the canonical Riemannian metric on $FM_i$ to obtain an
$O(n)$-invariant Riemannian metric on $FM_i$ to which we can apply
Proposition \ref{prop10}. From Proposition \ref{prop11},
there are \\
1. $O(n)$-basic
$\Z$-graded real vector bundles $\{ \check{E}_i \}_{i = 1}^\infty$ on
$\check{X}$,\\
2. Basic flat degree-$1$ superconnections $\{ \check{A}^\prime_i 
\}_{i = 1}^\infty$ on the $\check{E}_i$'s and
\\
3. $O(n)$-invariant Euclidean inner products $\{ h^{\check{E}_i} 
\}_{i = 1}^\infty$ on the $\check{E}_i$'s\\
so that for a given $j$ and large $i$, 
$\lambda_{p,j}(M_i)$ is $\epsilon$-close to 
$\lambda_{p,j}(X; E_i)$.  From
Proposition \ref{prop12}, after taking a subsequence we may assume that
the $\check{E}_i$'s are all isomorphic to a single 
$O(n)$-equivariant $\Z$-graded real
vector bundle $\check{E}$, the $E_{i,\alpha}$'s are all isomorphic to a
single $\Z$-graded real
vector bundle $E_\alpha$ and
the triples $\left\{ \left( {\frak I}_i,
\check{A}^\prime_i, h^{\check{E}_i} \right)
\right\}_{i=1}^\infty$ converge, after gauge transformations, to a triple
$\left({\frak I}, \check{A}^\prime, h^{\check{E}} \right)$. 

We claim that for $\check{x} \in \check{X}$, 
$\check{K}_{i,\check{x}}$ does not degenerate as $i \rightarrow \infty$.
To see this, note first that as the isotropy group
$H \subset O(n)$ acts freely and
affinely on the nilmanifold fiber
$\check{Z}_{i,\check{x}} \: = \: \Gamma_i \backslash N_i$ of $FM_i \rightarrow
\check{X}$, $H$ is virtually abelian and $H_0$, the connected component
of the identity in $H$, is a subgroup of the torus 
$C(\Gamma_i) \backslash C(N_i)$. 
As in the discussion before (\ref{upper}), 
$H_0$ acts isometrically on 
$\check{Z}_{i,\check{x}}$, with its orbit isometric to $H_0 \subset O(n)$.
Now $\check{E}_{i,\check{x}} \: \cong \: \Lambda^*({\frak n}_i^*)$.
As in
\cite[(6.6)]{Lott (2001)}, the Hermitian metric $h^{\check{E}_{i,\check{x}}}$
gives an orthogonal decomposition
${\frak n}_i^* \: = \: {\frak c(n_i)}^* \oplus ({\frak c(n_i)}^*)^\perp$.
Then for ${\frak x} \in {\frak h}$, the action of ${\frak I}_{\frak x}$ on 
$\check{E}_{i,\check{x}} \: \cong \: \Lambda^*({\frak c(n_i)}^*) \:
\widehat{\otimes} \: \Lambda^*(({\frak c(n_i)}^*)^\perp)$ is given by interior 
multiplication by ${\frak x}$ on the first factor. By passing to a subsequence,
we may assume that $\dim( {\frak c(n_i)} )$ is constant in $i$. Recalling that
$h^{\check{E}_{i,\check{x}}}$ comes from the $L^2$-inner product on the
affine-parallel forms on $\check{Z}_{i,\check{x}}$, 
it follows that
the actions of $H_0$ on $\{\check{E}_{i,\check{x}} \}_{i=1}^\infty$ are
related by isomorphisms with norms that are uniformly bounded above and below.
Then for all $i$, $\check{K}_{\check{x}} \: \cong \: \check{K}_{i,\check{x}}$.

Hereafter we think of all of the triples
$\left\{ \left( {\frak I}_i,
\check{A}^\prime_i, h^{\check{E}_i} \right)
\right\}_{i=1}^\infty$ as living on the same vector bundle $\check{E}$.
From the
convergence of the ${\frak I}_i$'s to ${\frak I}$, 
for large $i$ there are 
$O(n)$-equivariant automorphisms ${\cal A}_i$ of $\check{E}$ which converge
to the identity in the $C^\infty$-topology so that after
conjugating 
$\left( {\frak I}_i, \check{A}^\prime_i, h^{\check{E}_i} \right)$ 
with ${\cal A}_i$,
we may assume that $\check{K}_{\check{x}}$ as computed with
${\frak I}_i$ is the same as when computed with ${\frak I}$.
Then there is an isomorphism $T_i$ from the basic forms
defined using ${\frak I}$ to the basic forms defined using
${\frak I}_i$ which, using the notation of (\ref{to1}), is given on
$\check{X}_\alpha$ by
\begin{equation}
T_i \:  = \: \prod_j 
\left( 1 \: - \: e({\frak x}_j^*) \: {\frak I}_{i,{\frak x}_j} \right) \: 
\prod_j 
\left( 1 \: + \: e({\frak x}_j^*) \: {\frak I}_{{\frak x}_j} \right).
\end{equation}
Now $T_i^{-1} \: \check{A}^\prime_i \:
T_i$ is a flat degree-$1$ superconnection which is basic with
respect to ${\frak I}$.
For large $i$, Lemma \ref{lemma??} implies that
$\lambda_{p,j}(M_i)$ is $2 \epsilon$-close to $\lambda_{p,j}(X; E)$, the
$j$-th eigenvalue of the Laplacian $\triangle^E$ on
$\bigoplus_{a+b=p} \Omega^a(X; E^b)$.
In particular,
for all $j \ge 0$, $\lambda_{p,j}(X; E) \: \le \: e^{(J+2)\epsilon} \:
\Lambda$.
However, as a consequence of Lemma \ref{lemma??}, $\triangle^E_p$ has
a discrete spectrum, which is a contradiction.\\
{\bf 2.} 
If $A^1_{p,j,K}(M) \: = \: \infty$
then there is a sequence 
$\{g_i\}_{i = 1}^\infty$ of Riemannian
metrics on $M$ with $\diam(M,g_i) = 1$ and 
$\parallel R^{M}(g_i) \parallel_\infty \le K$ such that
$\lim_{i \rightarrow \infty} \lambda_{p,j}(M,g_i) \: = \: \infty$. 
As above, we can
find a new metric $g_i^\prime$ on $M$ which is $\epsilon$-close to $g_i$
satisfying 
$\parallel \nabla^k R^{M}(g_i^\prime) \parallel_\infty \: \le \:
A_k(n, \epsilon)$. 
Fix $\epsilon$ to be, say, $\frac{1}{2}$ and relabel the metrics 
$\{g_i^\prime\}_{i=1}^\infty$ to be $\{g_i\}_{i=1}^\infty$.
From \cite{Dodziuk (1982)} or \cite[Lemma 3]{Lott (2001)}, 
we again have
$\lim_{i \rightarrow \infty} \lambda_{p,j}(M, g_i) \: = \: \infty$.
As above, by taking a subsequence, we can assume that there is an
$O(n)$-manifold $\check{X}$ so
that $d_{GH}^{O(n)} (FM_i, \check{X}) \: \le \: \frac{1}{i}$. In particular,
$\dim(X) \: > \: 0$. Also as above, taking a further subsequence, 
we can assume that there are $\check{E}$, ${\frak I}$,
$\check{A}^\prime$ and $h^{\check{E}}$ so that if
$\lambda_{p,j}(X; E) \: < \: \infty$ then for large $i$,
$\lambda_{p,j}(M, g_i)$ is $\epsilon$-close to $\lambda_{p,j}(X; E)$. 
Suppose that $\dim(X) \: \ge \: p$. 
As $E^0$ is the trivial $\R$-bundle on $X$, there is an inclusion
$\Omega^p(X) \subset \bigoplus_{a + b = p} \Omega^a(X; E^b)$.
The Laplacian on $\Omega^p(X)$ is unbounded and so
$\lambda_{p,j}(X; E) < \infty$. This contradicts the assumption that
$\lim_{i \rightarrow \infty} \lambda_{p,j}(M,g_i) \: = \: \infty$.
Thus $\dim(X) \: < \: p$. Similarly, if $\Omega^a(X; E^b) \: \ne \:
0$ for any $a$ and $b$ satisfying $a \: + \: b \: = \: p$ then
$\triangle^E_p$ is unbounded and we get a contradiction.  Using the
construction of $E$ in terms of affine-parallel differential forms and
the finiteness statement of Proposition \ref{prop12}, the claim follows.\\
{\bf 3.} 
If $p \in \{0,1\}$ and $A^1_{n,p,j,K} \: = \: \infty$
then there is a sequence 
$\{M_i\}_{i = 1}^\infty$ of connected closed $n$-dimensional Riemannian
manifolds with $\diam(M_i) = 1$ and 
$\parallel R^{M_i} \parallel_\infty \le K$ such that
$\lim_{i \rightarrow \infty} \lambda_{p,j}(M_i) \: = \: \infty$.
As above, we can assume that there is an
$O(n)$-manifold $\check{X}$ so
that $d_{GH}^{O(n)} (FM_i, \check{X}) \: \le \: \frac{1}{i}$.
Also as above,
we can assume that there are $\check{E}$, ${\frak I}$,
$\check{A}^\prime$ and $h^{\check{E}}$ so that if
$\lambda_{p,j}(X; E) \: < \: \infty$ then for large $i$,
$\lambda_{p,j}(M_i)$ is $\epsilon$-close to $\lambda_{p,j}(X; E)$. 
As $E^0$ is the trivial $\R$-bundle on $X$ and $\dim(X) \: \ge \: 1$, 
there is an inclusion
$0 \neq \Omega^p(X) \subset \bigoplus_{a + b = p} \Omega^a(X; E^b)$.
The Laplacian on $\Omega^p(X)$ is unbounded and so
$\lambda_{p,j}(X; E) < \infty$. This contradicts the assumption that
$\lim_{i \rightarrow \infty} \lambda_{p,j}(M_i) \: = \: \infty$.
$\square$ \\ \\
{\bf Proof of Theorem \ref{tcor1} : }
Consider Theorem \ref{th5}.2 in the case $p = 2$. Then $\dim(X) = 1$. 
There are two cases.\\
1. $X = S^1$. As $M$ is an affine fiber bundle over $S^1$, it has the
claimed structure.\\
2. $X = [0,1]$, thought of as a singular space with two strata.
With respect to the singular fibration $q : M \rightarrow X$, put $Z \: = \:
q^{-1}(1/2)$ and $Z_i \: = \: q^{-1}(i-1)$. Then $M$ has the claimed
structure. We remark that the mapping cone is just defined up to
homeomorphism, which is why we use the word ``homeomorphic'' in
Theorem \ref{tcor1}.2.
 $\square$ \\ \\
{\bf Example 4 } 
\cite[Section 8]{Colbois-Ghanaat-Ruh (1999)} {\bf : } Let $N$ be the flat 
$3$-manifold ${\cal G}_6$ in the notation of
\cite[p. 122]{Wolf (1967)}. It has the rational homology of a
$3$-sphere. Let $h$ be a flat Riemannian metric on $N$.
Let $C > 0$ be the lowest
eigenvalue of $\triangle^N_1$.
Put $M \: = \: S^1 \times N$.
For $t \: > \: 0$, give $M$ the (flat) product metric 
$g_t \: = \: d\theta^2 \: + \: t^2 \: h$. Then 
one finds that the lowest eigenvalue of $\triangle^M_2$
is $C \: t^{-2}$. Taking $t \rightarrow 0$, we obtain that
$A^1_{2,j,K}(M) \: = \: \infty$ for all $j \in
\Z^+$ and $K \: \ge \: 0$. This is an example of Theorem \ref{tcor1}
in which $X \: = \: S^1$ and $Z \: = \: N$.

\section{Small Positive Eigenvalues} \label{sect8}

In this section we characterize the manifolds $M$ for which 
the $p$-form Laplacian has small positive eigenvalues.
We use the compactness result of
Section \ref{sect7} to show that if $M$ has $j$ small eigenvalues of
the $p$-form Laplacian, with $j > \bb_p(M)$, then $M$ collapses and
there is an associated basic flat degree-$1$ superconnection 
$\check{A}^\prime_\infty$
with $\dim(\HH^p(A^\prime_\infty)) \ge j$. We then use the spectral sequence of
$\check{A}^\prime_\infty$ to characterize when this can happen.
In Theorem \ref{tcor2}
we give a bound on the number of small eigenvalues
of the $1$-form Laplacian. 
In Theorems \ref{tcor4} and \ref{tcor7} we give bounds on the number of
small eigenvalues of the $p$-form Laplacian when one is sufficiently close
to a limit space of dimension $\dim(M)-1$ and characterize when
small eigenvalues can occur.

\begin{proposition} \label{prop14}
If $a_{p,j,K}(M) = 0$ and $j > \bb_p(M)$ then there are\\
1. An $O(n)$-equivariant affine fiber bundle
$FM \rightarrow \check{X}$,\\
2. A corresponding $O(n)$-equivariant $\Z$-graded real vector
bundle $\check{E} \rightarrow \check{X}$ and \\
3. A basic flat degree-$1$ superconnection 
$\check{A}^\prime_\infty$ on $\check{E}$\\
such that $\dim \left( \HH^p(A^\prime_\infty) \right) \: \ge \: j$.
\end{proposition}
\begin{pf}
Put $n = \dim(M)$.
If $a_{p,j,K}(M) = 0$ then there is a sequence $\{g_i\}_{i=1}^\infty$ in
${\cal M}(M, K)$ such that $\lim_{i \rightarrow \infty}
\lambda_{p,j}(M, g_i) = 0$. Let $M_i$ denote $M$ with the Riemannian metric
$g_i$. As $j > \bb_p(M)$, the $M_i$'s must collapse.
As in the proof of Theorem \ref{th5}, after smoothing the metrics and taking
a subsequence, we may assume that there is a smooth Riemannian $O(n)$-manifold
$\check{X}$ so that $\lim_{i \rightarrow \infty} d_{GH}^{O(n)}(FM_i, 
\check{X}) = 0$. From Proposition \ref{prop11}, there are \\
1. $O(n)$-basic
$\Z$-graded real vector bundles $\{ \check{E}_i \}_{i = 1}^\infty$ on
$\check{X}$,\\
2. Basic flat degree-$1$ superconnections $\{ \check{A}^\prime_i 
\}_{i = 1}^\infty$ on the $\check{E}_i$'s and
\\
3. $O(n)$-invariant Euclidean inner products $\{ h^{\check{E}_i} 
\}_{i = 1}^\infty$ on the $\check{E}_i$'s\\
so that $\lim_{i \rightarrow \infty} \lambda_{p,j}(X; E_i) = 0$.  From
Proposition \ref{prop12}, after taking a subsequence we may assume that
the $\check{E}_i$'s are all isomorphic to a single 
$O(n)$-equivariant $\Z$-graded real
vector bundle $\check{E}$, the $E_{i,\alpha}$'s are all isomorphic to a single
$\Z$-graded real vector bundle $E_\alpha$ and
the triples $\left\{ \left(
{\frak I}_i, \check{A}^\prime_i, h^{\check{E}_i} \right)
\right\}_{i=1}^\infty$ lie
in a compact subset $C$ of 
$({\cal S}_{\check{E},basic} \times 
{\cal H}_{\check{E},basic})/
{\cal G}_{\check{E},basic}$. 
Following the method of proof of Theorem \ref{th5}, we obtain a limit triple
$\left( {\frak I}_\infty, \check{A}^\prime_\infty,
h^{\check{E}_\infty} \right) \in {\cal S}_{\check{E},basic} \times 
{\cal H}_{\check{E},basic}$.

Let $g^{T\check{X}_\infty}$ denote the Riemannian metric on
$\check{X}$.
Let $\triangle^{E_\infty}$ denote the basic Laplacian constructed from 
$\check{A}^\prime_\infty$, $g^{T\check{X}_\infty}$ and $h^{\check{E}_\infty}$.
From Lemma \ref{lemma??} and the discreteness of the
spectrum of $\triangle^{M_i}$, it follows that the spectrum of
$\triangle^{E_\infty}$ is discrete. From the continuity of $\lambda_{p,j}$ as
a function on  ${\cal S}_{\check{E},basic} \times 
{\cal H}_{\check{E},basic}$,
we know that $\dim(\Ker(\triangle^{E_\infty}_p)) \:
\ge \: j$. It remains to show that $\Ker(\triangle^{E_\infty}_p) \cong
\HH^p(A^\prime_\infty)$. This amounts to a regularity issue.

In what follows, we omit the subscript $p$. 
Fix  $i$, with $\check{E}_i$ isomorphic as an $O(n)$-equivariant
vector bundle to
$\check{E}$. The idea will be to transfer the analysis of the superconnection
Laplacian to a more standard analysis on $M_i$.
Let $g^{T\check{X}_i}$ denote the Riemannian metric on
$\check{X}$ coming from the Riemannian affine fiber bundle
$FM_i \rightarrow \check{X}$ and let
$h^{\check{E}_i}$ denote the inner product on $\check{E}_i$ induced from the
Riemannian affine fiber bundle.

We know that 
$\Ker(\triangle^{E_\infty}) \: \cong \: \Ker(\check{A}^\prime_\infty)/
\overline{\Image(\check{A}^\prime_\infty)}$, where
$\check{A}^\prime_\infty$ acts on $\Omega_{basic, max} \left(
\check{X}; \check{E} \right)$. As in the proof of Theorem \ref{th5}.1,
we can conjugate $\left( {\frak I}_\infty, \check{A}^\prime_\infty, 
h^{\check{E}_\infty} \right)$
by an $O(n)$-equivariant automorphism of $\check{E}$ to make 
$\check{K}_{\check{x}}$ the same whether computed with 
${\frak I}_i$ or ${\frak I}_\infty$.
There is an isomorphism
from the basic forms defined using ${\frak I}_\infty$ to the basic
forms defined using ${\frak I}_i$ which, using the notation of
(\ref{to1}), is given on $\check{X}_\alpha$ by
\begin{equation}
T \:  = \: \prod_j 
\left( 1 \: - \: e({\frak x}_j^*) \: {\frak I}_{i,{\frak x}_j} \right) \: 
\prod_j 
\left( 1 \: + \: e({\frak x}_j^*) \: {\frak I}_{\infty,{\frak x}_j} \right).
\end{equation}
Then $T \: \check{A}^\prime_\infty \:
T^{-1}$ is a flat superconnection which is basic with
respect to ${\frak I}_i$.
Replacing $\check{A}^\prime_\infty$ and $h^{\check{E}_\infty}$ by
$T \: \check{A}^\prime_\infty \: T^{-1}$ and the inner product induced by $T$,
we may assume that the basic structure on $\check{E}$ is that of
${\frak I}_i$. We relabel $T \: \check{A}^\prime_\infty \: T^{-1}$ as
$\check{A}^\prime_\infty$.

As the underlying topological vector
space of the Hilbert space $\Omega_{basic, L^2} \left(
\check{X}; \check{E} \right)$ is the same whether the Hilbert space is
constructed using $\left( g^{T\check{X}_\infty}, h^{\check{E}_\infty} \right)$
or $\left( g^{T\check{X}_i}, h^{\check{E}_i} \right)$, it is equivalent to
consider the reduced cohomology of $\check{A}^\prime_\infty$ on 
$\Omega_{basic, max} \left(\check{X}; \check{E} \right)$, where the Hilbert 
space structure now
comes from $\left( g^{T\check{X}_i}, h^{\check{E}_i} \right)$.
Let $\triangle^E$ denote the basic Laplacian constructed from
$\check{A}^\prime_\infty$, ${\frak I}_i$,
$g^{T\check{X}_i}$ and $h^{\check{E}_i}$.

Let us first consider $\triangle^{E_i}$, the
basic Laplacian constructed
using $\check{A}^\prime_i$, ${\frak I}_i$, $g^{T\check{X}_i}$ and 
$h^{\check{E}_i}$. As in the proof of Proposition \ref{prop10},
there are isomorphisms
\begin{equation} \label{iso1}
\Omega^*(M_i) \: \cong \: \Omega^*_{basic}(FM_i) \: = \:
\Omega_{basic}(\check{X}; \check{E}_i) \: \oplus \: 
\Ker({\cal P}^{fib}).
\end{equation}
With respect to the inner products, we have isometric isomorphisms
\begin{equation} \label{withref}
\Omega^*_{L^2}(M_i, \mu_i) \: \cong \: \Omega^*_{basic, L^2}(FM_i) \: = \:
\Omega_{basic, L^2}(\check{X}; \check{E}_i) \: \oplus \: 
\Ker({\cal P}^{fib}),
\end{equation}
in terms of which the Laplacians are related by
\begin{equation}
\triangle^{M_i}_{\mu_i} \: \cong \: \triangle^{FM_i}_{basic} \: = \:
\triangle^{E_i} \: \oplus \: \triangle^{FM_i}_{basic} 
\big|_{Ker({\cal P}^{fib})}. 
\end{equation}
By standard elliptic theory,
$(I \: + \: \triangle^{M_i}_{\mu_i})^{-1}$ and
$(I \: + \: \triangle^{M_i}_{\mu_i})^{-1} \: d_{M_i}^*$ are compact,
and
$d_{M_i} \: (I \: + \: \triangle^{M_i}_{\mu_i})^{-1/2}$ is bounded.
Then $(I \: + \: \triangle^{E_i})^{-1}$ and
$(I \: + \: \triangle^{E_i})^{-1} \: 
\left( \check{A}^\prime_i \right)_{basic}^*$ are compact, and
$\check{A}^\prime_i \: (I \: + \: \triangle^{E_i})^{-1/2}$ is bounded.

Now $\triangle^{E}$ is also well-defined on 
$\Omega_{basic}(\check{X}; \check{E}_i)$. Hereafter, we will change notation
from $\check{E_i}$ to $\check{E}$. Put $y \: = \:
\check{A}^\prime_\infty \: - \: \check{A}^\prime_i \in
\Omega_{basic}(\check{X}; \End(\check{E}))$, which in particular is
a bounded operator on 
$\Omega_{basic, L^2} \left( \check{X}; \check{E} \right)$.
Then 
\begin{equation}
\triangle^E \: - \: \triangle^{E_i} \: = \:
\left( \check{A}^\prime_i \right)_{basic}^* y \: + \: y^* \check{A}^\prime_i
\: + \: y^* y.
\end{equation}
It follows that 
$(I \: + \: \triangle^{E_i})^{-1} \: \left( 
\triangle^E \: - \: \triangle^{E_i} \right) \:  
(I \: + \: \triangle^{E_i})^{-1/2}$ is compact. Then from 
\cite[Vol. IV, Chapter 13, Pf. of Corollary 4, p. 116]{Reed-Simon (1978)},
it follows that $\triangle^E$ has the same essential spectrum as
$\triangle^{E_i}$, i.e. the empty set, showing that 
$\triangle^E$ has a discrete spectrum.  Thus
$\check{A}^\prime_\infty$ has a closed image on
$\Omega_{basic, max} \left( \check{X}; \check{E} \right)$. Hence
\begin{equation}
\Ker \left( \check{A}^\prime_\infty \right)/
\overline{\Image \left( \check{A}^\prime_\infty \right)}
\: = \:
\Ker \left( \check{A}^\prime_\infty \right)/
\Image \left( \check{A}^\prime_\infty \right),
\end{equation}
where the latter is the (unreduced) cohomology of $\check{A}^\prime_\infty$ on
$\Omega_{basic, max} \left( \check{X}; \check{E} \right)$. 
(This also follows from the discreteness of the spectrum of
$\triangle^{E_\infty}$.) It remains
to show that this is isomorphic to the cohomology of 
$\check{A}^\prime_\infty$ on the smooth forms
$\Omega_{basic} \left( \check{X}; \check{E} \right)$.

There is an obvious cochain inclusion 
$\Omega_{basic} \left( \check{X}; \check{E} \right) \rightarrow
\Omega_{basic, max} \left( \check{X}; \check{E} \right)$.
We will construct a linear map $K$ on 
$\Omega_{basic, max} \left( \check{X}; \check{E} \right)$ of degree
$-1$ so that $I \: - \: \check{A}^\prime_\infty \: K \: - \:
K \: \check{A}^\prime_\infty$ sends 
$\Omega_{basic, max} \left( \check{X}; \check{E} \right)$ to
$\Omega_{basic} \left( \check{X}; \check{E} \right)$. This will give
a cochain homotopy equivalence between 
$\Omega_{basic} \left( \check{X}; \check{E} \right)$ and
$\Omega_{basic, max} \left( \check{X}; \check{E} \right)$, showing that
the two complexes have isomorphic cohomologies.

To the affine fiber bundle $FM_i \rightarrow \check{X}$ we associate
an infinite-dimensional 
$\Z$-graded $O(n)$-basic vector bundle $\check{W}$ on $\check{X}$, as in 
\cite[Section 5]{Lott (2001)}, so that 
$\Omega \left( \check{X}; \check{W} \right) \cong \Omega^*(FM_i)$.
The inclusion of fibers $\check{E}_{\check{x}} \subset \check{W}_{\check{x}}$
is isomorphic to the inclusion $1 \otimes \Lambda^*({\frak n}^*) \subset 
C^\infty(\check{Z}_{\check{x}}) 
\otimes \Lambda^*({\frak n}^*)$. Then the inclusion
$\Id \: \otimes \: \End(\Lambda^*({\frak n}^*)) \subset 
\End( C^\infty(\check{Z}_{\check{x}}) \otimes \Lambda^*({\frak n}^*))$ induces
an extension
$Y \in \Omega \left( \check{X}; \End (\check{W}) \right)$ of $y$.
As $y \in \Omega_{basic} \left( \check{X}; \End (\check{E}) \right)$, it
follows that
$Y \in \Omega_{basic} \left( \check{X}; \End (\check{W}) \right)$ and
that the corresponding map on 
$\Omega_{basic} \left( \check{X}; \check{W} \right) \cong \Omega^*(M_i)$
is diagonal with respect to (\ref{iso1}), so that we can write
$Y \: = \: 
\begin{pmatrix}
y & 0 \\
0 & *
\end{pmatrix}$. 

Put $D \: = \: d_{M_i} \: + \: Y$, a pseudodifferential operator on 
$M_i$ of order $1$. It decomposes with respect to (\ref{withref}) as
$D \: = \: 
\begin{pmatrix}
\check{A}^\prime_\infty & 0 \\
0 & *
\end{pmatrix}$. Then $D D^* \: + \: D^* D$ is a elliptic
pseudodifferential operator on $M_i$ of order $2$, which decomposes as
$D D^* \: + \: D^* D\: = \: 
\begin{pmatrix}
\triangle^E & 0 \\
0 & *
\end{pmatrix}$. Fix $t \: > \: 0$ and
let $K$ be the restriction of $D^* \: 
\frac{I \: - \: e^{- \: t(DD^* \: + \: D^*D)}}{DD^* \: + \: D^* D}$ 
to the first factor
$\Omega_{basic, L^2}(\check{X}; \check{E})$. 
Then
$K \: = \: (\check{A}^\prime_\infty)^*_{basic} \frac{I \: - \: 
e^{- \: t \triangle^E}}{\triangle^E}$, as defined spectrally.
We have
\begin{equation}
I \: - \: \check{A}^\prime_\infty \: K \: - \: K \:
\check{A}^\prime_\infty \: = \: e^{- \: t \triangle^E}.
\end{equation}
Now $e^{- \: t \triangle^E}$
is the restriction of $e^{- \: t (DD^* \: + \: D^*D)}$ to
$\Omega_{basic, L^2}(\check{X}; \check{E})$.
By elliptic theory,
$e^{- \: t(DD^* \: + \: D^*D)}$ maps 
$\Omega^*_{L^2}(M_i, \mu_i)$ to $\Omega^*(M_i)$. Hence 
$e^{- \: t \triangle^E}$
maps $\Omega_{basic, L^2}(\check{X}; \check{E})$ to
\begin{equation}
\Omega_{basic, L^2}(\check{X}; \check{E}) \: \cap \: \Omega^*(M_i) \: = \:
\Omega_{basic}(\check{X}; \check{E}).
\end{equation}
This proves the proposition. 
\end{pf}

Under the hypotheses of Proposition
\ref{prop14}, using the spectral sequence of Section \ref{sect6} we deduce that
\begin{equation} \label{eq8.8}
j \: \le \: \dim( \HH^p(A^\prime_\infty)) \: \le \:
\sum_{a+b=p} \dim \left( \HH^a(X; \HH^b(A_{\infty,[0]}^\prime)) \right). 
\end{equation}

Given an $O(n)$-equivariant affine
fiber bundle $FM \rightarrow \check{X}$ with fiber
$\check{Z} \: = \: \Gamma \backslash N$, there is a spectral
sequence to compute $\HH^*(M; \R)$ as the cohomology of a basic flat 
degree-$1$ superconnection
$\check{A}^\prime$ as in \cite[Example 4]{Lott (2001)}, working
equivariantly with respect to the $O(n)$-action and using basic forms.
With the notation of Section \ref{sect6}, the 
$E_2$-term of this spectral sequence is
$E_2^{p,q} \: = \: \HH^p(X; \HH^q(A^\prime_{[0]}))$.
If $x \in X$ is covered by $\check{x} \in \check{X}$ then 
the stalk of $\HH^q(A^\prime_{[0]})$ at $x$ is given by the cohomology of
$\check{A}^\prime_{[0]} \: = \: d^{\check{Z}}$ on $\check{K}_{\check{x}}^H$.
If $\check{x} \in \check{X}$ has isotropy group $H \subset O(n)$ 
then for ${\frak x} \in
{\frak h}$, ${\frak I}_{\frak x}$ is interior multiplication by the
corresponding vector field on $\Omega^*(\check{Z}_{\check{x}})$. We see that
$\check{K}_{\check{x}}^H$ consists of the $H$-basic forms on the fiber
$\check{Z}_{\check{x}}$. Then the stalk of $\HH^q(A^\prime_{[0]})$ at $x$
is isomorphic to
$\HH^q(\check{Z}_{\check{x}}/H; \R)$; recall that $H$ acts freely on
$\check{Z}_{\check{x}}$.
 
On the other hand, there
is a spectral sequence to compute $\HH^*(M; \R)$ from the map 
$r : M \rightarrow X$ \cite[p. 179]{Bott-Tu (1982)}, with
$E_2$-term $\HH^*(X; \HH^*(Z; \R))$. Here $\HH^*(Z; \R)$ is a sheaf on
$X$ whose stalk over $x \in X$ is $\HH^*(r^{-1}(x); \R)$.
As $r^{-1}(x) \: = \: \check{Z}_{\check{x}}/H$, we see that the
two spectral sequences have similar $E_2$-terms. In fact the
two spectral sequences are equivalent, as follows from the construction in
\cite[p. 179]{Bott-Tu (1982)}.

To obtain results about small eigenvalues, the idea now is to
compare the spectral sequence of $\check{A}^\prime_\infty$ with the
spectral sequence of the map $M \rightarrow X$.
 \\ \\
{\bf Example 5 : } As an illustration of the methods, we analyze the
behavior of the differential form Laplacian under the
collapse of $S^3$ to an interval which is described in
\cite[Example 1.5]{Cheeger-Gromov (1986)}.
With respect to the isometric action of 
$\SO(2) \times \SO(2) \: \subset \:
\SO(4)$ on the round
$S^3$, we shrink the metric on $S^3$ in the direction of a subgroup
$\R \subset \SO(2) \times \SO(2)$ 
of irrational slope. The resulting metrics approach the closed interval
$X \: = \: S^3/(\SO(2) \times \SO(2))$ 
in the Gromov-Hausdorff topology.

For simplicity, we consider the principal spin bundle $S^3 \times S^3$,
with structure group $G \: = \: \SU(2)$,
instead of the orthonormal frame bundle $FM$. Then during the collapse we have
$\check{X} \: = \: \C P^2 \: \# \: \overline{\C P^2}$ and 
$\check{Z} \: = \: T^2$. Correspondingly, $\check{E}$ is a vector bundle
over $\check{X}$ with fiber isomorphic to $\Lambda^*((\R^2)^*)$, on which
$\check{A}^\prime_{[0]}$ acts as the zero map. 
If $\check{x} \in \check{X}$ covers $x \in \Int(X)$ then its isotropy group $H$
is trivial, while if $\check{x} \in \check{X}$ covers $x \in \partial X$
then its isotropy group $H$ is isomorphic to $\U(1)$.
The sheaf $\HH^0( A^\prime_{\infty,[0]})$ 
has stalk $\R$ over each $x \in X$. The sheaf
$\HH^1( A^\prime_{\infty,[0]})$ has stalk $\R$ over $x \in \partial X$ and 
stalk
$\R^2$ over $x \in \Int(X)$, where we can think of the $\R$ over one component
of $\partial X$ as 
corresponding to $\R \oplus 0$ and the $\R$ over the other component of
$\partial X$ as corresponding to
$0 \oplus \R$. The sheaf
$\HH^2( A^\prime_{\infty,[0]})$ has stalk $0$ over $x \in \partial X$ and stalk
$\R$ over $x \in \Int(X)$. One finds that the only nonzero components of the
$E_2$-term of the spectral sequence for $\check{A}^\prime_{\infty}$ are
$\HH^0 \left( X ; \HH^0( A^\prime_{\infty,[0]}) \right) \: = \: 
\HH^1 \left( X ; \HH^2( A^\prime_{\infty,[0]}) \right) \: = \: \R$.
(These correspond to the zero-eigenvalues of the differential form
Laplacian on $S^3$.)

Let $r \: : \: S^3 \rightarrow X$ be the quotient map.
If $x \in \Int(X)$ then $r^{-1}(x) \: = \: T^2$, while
if $x \in \partial X$ then $r^{-1}(x) \: = \: S^1$. One finds that the
Leray spectral sequence to compute $\HH^*(S^3; \R)$ from the map $r$
coincides with the spectral sequence to compute the cohomology of
$\check{A}^\prime_{\infty}$.
The conclusion
is that there are no small positive eigenvalues in the collapse. This
is in contrast to what happens in the Berger collapse of $S^3$ to $S^2$
\cite[Example 1.2]{Colbois-Courtois (1990)}.
\begin{corollary} \label{cor8}
For any $K \ge 0$, $a_{0,2,K}(M) > 0$. That is, there are no small
positive eigenvalues of the Laplacian on functions.
\end{corollary}
\begin{pf}
Suppose that $a_{0,2,K}(M) = 0$. Using Proposition \ref{prop14} and
equation (\ref{eq8.8}) in the case $p = 0$, we conclude that
$2 \: \le \: \dim \left( \HH^0(X; \HH^0(A_{\infty,[0]}^\prime)) \right)$. 
However, $E^0$ is the trivial $\R$-bundle on $X$, on which
$A^\prime_{\infty,[0]}$ acts by zero. Then
$\HH^0(A^\prime_{\infty, [0]})$ is the trivial $\R$-bundle on $X$ and so
$\dim(\HH^0(X; \HH^0(A^\prime_{\infty,[0]}))) = 1$, which is a contradiction.
\end{pf}
{\bf Remark :} Corollary \ref{cor8} is true under the weaker
assumption of a lower bound on the Ricci curvature (\cite{Berard (1988)} and
references therein).
\begin{corollary} \label{tcor3}
Suppose that
$a_{p,j,K}(M) = 0$ and $j > \bb_p(M)$.  Let $X$ be the limit space
coming from the above argument.
If $\dim(X) = 0$, write the almost flat
manifold $M$ topologically
as the quotient of a nilmanifold $\widehat{\Gamma} \backslash N$
by a finite group $F$. Let ${\frak n}$ denote the Lie algebra of the nilpotent
Lie group $N$ and let $\cdot^F$ denote $F$-invariants.
Then $j \: \le \: \dim(\Lambda^p({\frak n}^*)^F)$.
\end{corollary}
\begin{pf}
With reference to Proposition \ref{prop14},
in this case $X$ is a point and $E = \Lambda^*({\frak n}^*)^F$. Then for
any superconnection $A^\prime_\infty$ on $E$, i.e. differential 
$A^\prime_{\infty,[0]}$ on $E$, we have 
$\dim(\HH^p(A^\prime_\infty)) \: \le \: 
\dim(E^p)$.
\end{pf}
{\bf Remark : }
It follows from \cite[Theorem 6]{Lott (2001)}
that in fact $a_{p,j,K}(M) = 0$ for 
$j \: = \: \dim(\Lambda^p({\frak n}^*)^F)$. \\ \\
{\bf Proof of Theorem \ref{tcor2} : }
With reference to Proposition \ref{prop14} and equation (\ref{eq8.8}),
\begin{equation}
j \: \le \: \dim \left( \HH^1(X; \HH^0(A^\prime_{\infty,[0]}) \right) +
\dim \left( \HH^0(X; \HH^1(A^\prime_{\infty,[0]}) \right).
\end{equation}
As $\HH^0(A^\prime_{\infty,[0]})$ is the trivial $\R$-bundle on $X$,
$\dim \left( \HH^1(X; \HH^0(A^\prime_{\infty,[0]}) \right) = \bb_1(X)$.
As $A^\prime_{\infty,[0]}$ acts by zero on $E^0$, 
there is an injection $\HH^1(A^\prime_{\infty,[0]}) \rightarrow E^1$. 
Let $\dim \left( \HH^1(A^\prime_{\infty,[0]}) \right)$ denote the
dimension of the generic stalk of the sheaf 
$\HH^1(A^\prime_{\infty,[0]})$ and put
$\dim(E^1) \: = \: \dim(E^1_{\beta})$, where $\beta$ again denotes the
principal normal orbit type.
Then 
\begin{equation} \label{eq8.10}
\dim \left( \HH^0(X; \HH^1(A^\prime_{\infty,[0]})) \right) \: \le \:
\dim \left( \HH^1(A^\prime_{\infty,[0]}) \right) 
\le \: \dim \left( E^1 \right)
\: \le \: \dim(M) \: - \: \dim(X).
\end{equation}
Thus $j \: \le \: \bb_1(X) \: + \dim(M) \: - \: \dim(X)$.
On the other hand, the spectral sequence for $\HH^*(M; \R)$ gives
\begin{equation} \label{eq8.11}
\HH^1(M; \R) \: = \: \HH^1(X; \R) \: \oplus \: \Ker \left( 
\HH^0(X; \HH^1(Z; \R)) \rightarrow \HH^2(X; \R) \right).
\end{equation}
In particular, $\bb_1(X) \: \le \: \bb_1(M)$.
The corollary follows. $\square$ \\ \\
{\bf Remark : } Using heat equation methods \cite{Berard (1988)} one can
show that there is an increasing function $f$ such that if $\Ric(M) \: \ge \:
- (n-1) \: \lambda^2$ and $\diam(M) \: \le \: D$ then the number of small 
eigenvalues of the $1$-form Laplacian is bounded above by 
$f(\lambda \: D)$. This result is weaker than Theorem \ref{tcor2} when
applied to manifolds with sectional curvature bounds, but is more general
in that it applies to manifolds with just a lower Ricci curvature bound.\\ \\
{\bf Proof of Theorem \ref{tcor4} :}
In general, if  $X$ is a compact metric space of Hausdorff dimension $n-1$, and
a manifold $M^n$ with $\parallel R^M \parallel_\infty \: \le \: K$ 
is sufficiently 
Gromov-Hausdorff close to $X$, then 
the description of the local geometry of $M$ in 
\cite[Theorem 1.3]{Cheeger-Fukaya-Gromov (1992)} says that $M$ has the
following structure. 
Given $m \in M$, there are a neighborhood $U_m$ of $m$, a finite
regular covering $\widehat{U}_m$ of $U_m$ with covering group $F$ and
a locally free 
$S^1$-action on $\widehat{U}_m$ which is $F$-equivariant with respect
to a homomorphism $\eta : F \rightarrow \Aut(S^1) (\cong \Z_2)$. 
(In fact, we can take $F$ to be $\{e\}$ or $\Z_2$.) Then
we can take $X$ so that locally, it is
the quotient $\widehat{U}_m/(F \widetilde{\times}
S^1)$.  Hence $X$ is an orbifold and $M \rightarrow X$ is an orbifold
circle bundle, with its orientation bundle ${\cal O}$ given locally as
$\widehat{U}_m \times_{F \widetilde{\times} S^1} \R \rightarrow
\widehat{U}_m / (F \widetilde{\times} S^1)$. Here
$F \widetilde{\times} S^1$ acts on $\R$ through $F$.

Suppose that the claim of the corollary is not true.
Then there is a sequence of connected closed $n$-dimensional Riemannian 
manifolds $\{(M_i, g_i)\}_{i=1}^\infty$ with $\parallel R^{M_i}(g_i) 
\parallel_\infty \: \le \: K$ and $\lim_{i \rightarrow \infty} M_i = X$
which provides a counterexample.
As there is a finite number of isomorphism
classes of flat real line
bundles on $X$, after passing to a subsequence we may assume that
each $M_i$ is an orbifold circle bundle over $X$ with a fixed 
orientation bundle
${\cal O}$ and that $\lim_{i \rightarrow \infty} \lambda_{p,j}(M_i, g_i) = 0$
for $j \: = \: \bb_p(X) \: + \: \bb_{p-1}(X; {\cal O}) \: + \: 1$. 
As in the proofs of Theorem \ref{th5} and Proposition \ref{prop14}, 
we obtain $E \: = \: E^0 \: \oplus \: E^1$ on $X$,
with $E^0$ a trivial $\R$-bundle and $E^1 \: = \: {\cal O}$, and a limit
superconnection $A^\prime_\infty$ on $E$ with $A^\prime_{\infty,[0]} = 0$ and
$A^\prime_{\infty,[1]} = \nabla^E$, the canonical flat connection. Then as in
(\ref{eq8.8}), we obtain
\begin{equation}
j \: \le \: \dim( \HH^p(A^\prime_\infty)) \: \le \:
\bb_p(X) \: + \: \bb_{p-1}(X; {\cal O}),
\end{equation}
which is a contradiction. \\ \\
{\bf Proof of Theorem \ref{th6} :} The proof of the theorem is along the
same lines as that of \cite[Theorem 5]{Lott (2001)}, replacing
\cite[(7.8)]{Lott (2001)} by (\ref{eq8.8}). We omit the details.\\ \\
{\bf Proof of Theorem \ref{tcor7} :}
The $E_2$-term of the spectral sequence to compute $\HH^*(M; \R)$ consists
of $E_2^{p,0} = \HH^p(X; \R)$ and $E_2^{p,1} = \HH^p(X; {\cal O})$. The
differential is ${\cal M}_{\chi}$. The
corollary now follows from Theorem \ref{th6}. $\square$

\end{document}